\documentclass[smallextended,envcountsect]{svjour3} 
\smartqed 
\usepackage{graphicx}
\usepackage[english]{babel} 
\usepackage[utf8]{inputenc}
\usepackage{amsmath} 
\usepackage{amssymb}
\usepackage{mathrsfs}
\usepackage{amsfonts} 
\usepackage{stmaryrd}
\usepackage{marvosym}

\usepackage{algorithm}
\usepackage{algpseudocode}

\usepackage{graphicx}
\usepackage{psfrag}

\journalname{JOTA}

\newcommand{\T}{\mathscr{T}}
\newcommand{\Sides}{\mathscr{S}}

\usepackage[usenames,dvipsnames]{xcolor}
\newcommand{\FF}[1]{{\color{black}#1}}
\newcommand{\EO}[1]{{\color{black}#1}}

\begin{document}

\title{A Posteriori Error Estimates for an Optimal Control Problem with a Bilinear State Equation}

\author{Francisco Fuica$^{1}$  \and  Enrique Ot\'arola$^{1}$}

\institute{\Letter ~ Enrique Ot\'arola   \\
 enrique.otarola@usm.cl  \at
 			$^1$ Departamento de Matem\'atica, Universidad T\'ecnica Federico Santa Mar\'ia, Av. Espa\~{n}a 1680, Valpara\'iso, Chile.
}

\date{Received: date / Accepted: date}

\maketitle

\begin{abstract}
We propose and analyze a posteriori error estimators for an optimal control problem that involves an elliptic partial differential equation as state equation and a control variable that enters the state equation as a coefficient; \FF{pointwise constraints on the control variable are considered as well}. We consider two different strategies \FF{to approximate optimal variables}: a fully discrete scheme \FF{in which} the admissible control set is discretized with piecewise constant functions and a semi-discrete scheme where the admissible control set is not discretized; \FF{the latter scheme being based on the so-called variational discretization approach}. We design, for each solution technique, an a posteriori error estimator and show, in two and three dimensional Lipschitz \FF{polygonal/polyhedral} domains \FF{(not necessarily convex)}, that the proposed error estimator is reliable and efficient. \FF{We design, based on the devised estimators,} adaptive strategies that \FF{deliver} optimal experimental rates of convergence for the \FF{performed} numerical examples.
\end{abstract}

\keywords{optimal control problems \and bilinear equations \and finite elements \and a posteriori error estimates \and adaptive finite element methods}

\subclass{49M25	   % Discrete approximations
\and 65N15         % Error bounds
\and 65N30         % Finite elements, Rayleigh-Ritz and Galerkin methods,  finite methods
\and 65N50         % Mesh generation and refinement 
}

\section{Introduction}\label{sec:intro}

The \FF{development and study of discretization techniques, based on finite elements,} for distributed control--constrained linear--quadratic elliptic optimal control problems have been widely studied in the literature; see \cite{MR2843956,MR2516528} for an extensive list of references. These discretization techniques are mainly divided \FF{into} two categories, which rely on the discretization of the state and adjoint equations; they differ on whether or not the admissible control set is also discretized.
In contrast to these advances, the study of solution techniques for optimal control problems where the control variable enters the state equation as a coefficient is not as developed. One of the main sources of difficulty within these type of problems is that the solution of the state equation depends nonlinearly on the control variable  \cite{MR2536007}. Consequently, uniqueness of solutions cannot be guaranteed.

In this work, we will \FF{focus on the development} and analysis of efficient solution techniques for the optimal control problem \eqref{eq:minimize_cost_func}--\eqref{def:box_constraints}, which incorporates the control variable as a coefficient in the state equation; the control variable is not a source term. We immediately mention that this problem can be interpreted as a particular instance of parameter estimation. To the best of our knowledge, the first work that provides an analysis for suitable finite element discretizations for problem \eqref{eq:minimize_cost_func}--\eqref{def:box_constraints} is \cite{MR2536007}. In this work, the authors propose, on convex \FF{polygonal/polyhedral} domains and quasi--uniform meshes, two fully discrete schemes that discretize the admissible control set with piecewise constant and piecewise linear functions; the state and adjoint equations are discretized with piecewise linear functions. Estimates for the error committed within the approximation of a control variable are derived in \cite[Corollaries 5.6 and 5.10]{MR2536007}. In addition, an error estimate for a post-processing strategy is obtained in \cite[Theorem 5.18]{MR2536007}. The results obtained in \cite{MR2536007} were later extended to mixed and stabilized finite element methods in \cite{MR3103238} and \cite{MR3693332}, respectively.

\FF{A particular class of numerical methods that has proven a competitive performance when are used to approximate} solutions to PDE--constrained optimization problems, and the ones we shall \FF{consider} in this work, are adaptive finite element methods (AFEMs). AFEMs are iterative methods \FF{recognizable by their capability to improve} the quality of a \FF{discrete} approximation to a corresponding PDE while \FF{keeping} an \FF{efficient} distribution of computational resources. A \FF{crucial component} of an AFEM is an a posteriori error estimator, which is a computable quantity, \FF{depending} on the \FF{problem data} and discrete solution, \FF{that provides local} information about the quality of the approximate solution. The a posteriori error analysis for control--constrained linear--quadratic optimal control problems has achieved several advances in recent years. We refer the interested reader to \cite{MR1887737,MR1780911,MR2434065,MR3212590,MR3621827,MR4122501} for a discussion. As opposed to these advances, the analysis of AFEMs for optimal control problems involving nonlinear or bilinear equations is rather scarce. To the best of our knowledge, the work \cite{MR2680928} appears to be the first that provides a posteriori error estimates for  \eqref{eq:minimize_cost_func}--\eqref{def:box_constraints}. In this work, the authors develop a posteriori error estimators for two fully discrete approximation schemes of \eqref{eq:minimize_cost_func}--\eqref{def:box_constraints} and obtain global reliability estimates \cite[Theorems 4.1 and 4.3]{MR2680928} and \emph{global} efficiency results \cite[Lemmas 4.2 and 4.3]{MR2680928}. We also mention the work \cite{MR3174031}, where a posteriori error estimates for a parabolic version of \eqref{eq:minimize_cost_func}--\eqref{def:box_constraints} have been analyzed; an efficiency analysis, however, was not provided. \FF{We conclude this paragraph by mentioning the work \cite{MR2373479}, where the authors provide, on the basis of a posteriori error estimators, upper bounds for discretization errors with respect to a cost functional and with respect to a given quantity of interest; the latter being an arbitrary functional depending on the control and the state variables. In our work, we derive upper and lower bounds for the approximation error when is measured in an \emph{energy} norm; see below for a discussion.}

\FF{In the present manuscript, we consider} two different strategies to discretize the optimal control problem \eqref{eq:minimize_cost_func}--\eqref{def:box_constraints}: a semi-discrete scheme, \FF{based on the so-called variational discretization approach \cite{MR2122182}, in which} the admissible control set is not discretized, and a fully discrete scheme, where control variables are \FF{approximated} by using piecewise constant functions. \FF{We devise,} for each one of \FF{the aforementioned} schemes, a residual--based a posteriori error estimator. For the fully discrete scheme the error estimator is \FF{formed by} the sum of three contributions\FF{: two of them} are related to the discretization of the state and adjoint equations \FF{while the remaining one is related to the discretization of the} admissible control set. In contrast, the error estimator for the variational discretization approach is \FF{formed by} only two contributions that are related to the discretization of the state and adjoint equations. In two and three dimensional Lipschitz \FF{polygonal/polyhedral domains (not necessarily convex)}, we obtain reliability and efficiency estimates. 

In what follows we list what, we believe, are the main contributions of our work:
\begin{itemize}
\item For the fully and semi-discrete schemes that we consider, we devise a posteriori error estimators; both being different from the ones in \cite{MR2680928} \FF{and \cite{MR2373479}}.

\item For the aforementioned solution techniques, we prove that the corresponding local error indicators associated to the discretization of the state and adjoint equations are \emph{locally efficient}. This analysis improves the global one in \cite[Lemmas 4.2 and 4.3]{MR2680928}. We also prove that the total error \EO{indicator} associated to the variational discretization approach is \emph{locally efficient} \EO{(cf.~Theorem \ref{thm:global_eff_var})}; the one associated to the fully discrete scheme, as customary, being globally efficient \EO{(cf.~Theorem \ref{thm:global_eff})}. 

\item We design a simple adaptive loop that delivers optimal experimental rates of convergence for all the involved individual contributions of the corresponding error. The loop based on the a posteriori error indicators devised for the variational discretization approach delivers quadratic rates of convergence for the error approximation of a control variable. This substantially improves the approximation properties that can be achieved by the considered fully discrete scheme. The indicators devised for the latter scheme tend to refine the involved meshes in regions where the restrictions of the control variable become active. These DOFs seem not necessary for an accurate approximation of  control variables. An scheme based on piecewise linear approximation of the admissible control set would suffer the same limitations in terms of avoidable refinement \cite{MR4122501}.

\end{itemize}

The rest of the paper is organized as follows. In section \ref{sec:the_problem_and_notation} we introduce the optimal control problem under consideration and set notation. Basic results for the state equation as well as basic a posteriori error estimates are reviewed in section \ref{sec:review_bilinear}.
In section \ref{sec:the_ocp} we review the existence of solutions for the optimal control problem as well as first and second order optimality conditions. The crucial part of our work are sections \ref{sec:a_posteriori_fully} and \ref{sec:a_posteriori_semi}, where we design and analyze a posteriori error estimators for the fully and semi-discrete schemes, respectively. Finally, in section \ref{sec:numerical_ex} we present numerical examples in two and three dimensional domains that illustrate the theory and reveal a competitive performance of the devised AFEMs.

\section{The Problem and Notation}\label{sec:the_problem_and_notation}

Let us precisely introduce the optimal control problem that will be considered in our work, set notation, and describe the setting we shall operate with.

\subsection{Presentation of the Problem}

In this work we are interested in the design and analysis of a posteriori error estimates for an optimal control problem governed by an elliptic partial differential equation (PDE) as state equation. Our main source of difficulty here is that the control variable enters the state equation as a coefficient; control constraints are also considered. Let us make this discussion precise. Let $\Omega\subset\mathbb{R}^d$, with $d\in \{2, 3\}$, be an open and bounded \FF{polygonal/polyhedral} domain with Lipschitz boundary $\partial\Omega$ \FF{\cite[Chapter 4]{MR2424078}}. Given a desired state $y_\Omega\in L^2(\Omega)$ and a regularization parameter $\alpha>0$, let us introduce the cost functional
\begin{equation}\label{def:cost_functional}
J(y,u):=\frac{1}{2}\|y-y^{}_\Omega\|_{L^2(\Omega)}^2+\frac{\alpha}{2}\|u\|_{L^2(\Omega)}^2.
\end{equation}
We are thus interested in the following optimal control problem: Find
\begin{equation}\label{eq:minimize_cost_func}
\min{J(y,u)}
\end{equation}
subject to the elliptic PDE 
\begin{equation}\label{eq:state_equation}
-\Delta y+uy=f \text{ in } \Omega,\qquad y= 0  \text{ on } \partial\Omega,
\end{equation}
where $f\in L^2(\Omega)$ denotes an external source, and the control constraints 
\begin{equation}\label{def:box_constraints}
u\in\mathbb{U}_{ad}, \qquad \mathbb{U}_{ad}:=\{ v\in L^2(\Omega): 0<\texttt{a} \leq v \leq \texttt{b} \text{ a.e. in } \Omega \}.
\end{equation}
Here, $\texttt{a},\texttt{b} \in  \mathbb{R}^{+}$ satisfy $\texttt{a}<\texttt{b}$. 

%%%%%%%%%%%%%%%%%%%%%%%%%%%%%%%%%%%%%%%%%%%%%%%%%%%%%%%%%%%%
%%%%%%%%%%%%%%%%%%%%%%%%%%%%%%%%%%%%%%%%%%%%%%%%%%%%%%%%%%%%
%%%%%%%%%%%%%%%%%%%%%%%%%%%%%%%%%%%%%%%%%%%%%%%%%%%%%%%%%%%%
%%%%%%%%%%%%%%%%%%%%%%%%%%%%%%%%%%%%%%%%%%%%%%%%%%%%%%%%%%%%

\subsection{Notation}\label{sec:notation}

Let us set notation and describe the setting we shall operate with. Throughout this work $d\in\{2,3\}$ and $\Omega\subset\mathbb{R}^d$ is an open and bounded \FF{polygonal/polyhedral} domain with Lipschitz boundary $\partial\Omega$ \FF{\cite[Chapter 4]{MR2424078}}. Notice that we do not assume that $\Omega$ is convex. If $\mathcal{X}$ and $\mathcal{Y}$ are normed vector spaces, we write $\mathcal{X}\hookrightarrow\mathcal{Y}$ to denote that $\mathcal{X}$ is continuously embedded in $\mathcal{Y}$. We denote by $\|\cdot\|_\mathcal{X}$ the norm of $\mathcal{X}$. The relation $a \lesssim b$ indicates that $a \leq C b$, with a positive constant that depends neither on $a$, $b$ nor on the involved discretization parameters. The value of $C$ might change at each occurrence.

%%%%%%%%%%%%%%%%%%%%%%%%%%%%%%%%%%%%%%%%%%%%%%%%%%%%%%%%%%%%
%%%%%%%%%%%%%%%%%%%%%%%%%%%%%%%%%%%%%%%%%%%%%%%%%%%%%%%%%%%%
%%%%%%%%%%%%%%%%%%%%%%%%%%%%%%%%%%%%%%%%%%%%%%%%%%%%%%%%%%%%
%%%%%%%%%%%%%%%%%%%%%%%%%%%%%%%%%%%%%%%%%%%%%%%%%%%%%%%%%%%%
%%%%%%%%%%%%%%%%%%%%%%%%%%%%%%%%%%%%%%%%%%%%%%%%%%%%%%%%%%%%
%%%%%%%%%%%%%%%%%%%%%%%%%%%%%%%%%%%%%%%%%%%%%%%%%%%%%%%%%%%%
%%%%%%%%%%%%%%%%%%%%%%%%%%%%%%%%%%%%%%%%%%%%%%%%%%%%%%%%%%%%
%%%%%%%%%%%%%%%%%%%%%%%%%%%%%%%%%%%%%%%%%%%%%%%%%%%%%%%%%%%%

\section{The State Equation}
\label{sec:review_bilinear}

In this section, we briefly review some results related to the well-posedness of problem \eqref{eq:state_equation}. Additionally, we \FF{present} a posteriori error estimates for a \FF{specific} finite element setting.

%%%%%%%%%%%%%%%%%%%%%%%%%%%%%%%%%%%%%%%%%%%%%%%%%%%%%%%%%%%%
%%%%%%%%%%%%%%%%%%%%%%%%%%%%%%%%%%%%%%%%%%%%%%%%%%%%%%%%%%%%
%%%%%%%%%%%%%%%%%%%%%%%%%%%%%%%%%%%%%%%%%%%%%%%%%%%%%%%%%%%%
%%%%%%%%%%%%%%%%%%%%%%%%%%%%%%%%%%%%%%%%%%%%%%%%%%%%%%%%%%%%

\subsection{Weak Formulation}

Let $\mathfrak{f}$ be a given forcing term in $L^2(\Omega)$ and $\mathfrak{u}$ be an arbitrary function in $\mathbb{U}_{ad}$. With this setting at hand, we introduce the following weak problem:
\begin{equation}\label{eq:weak_eq}
z\in H_0^1(\Omega): 
\quad
(\nabla z,\nabla v)_{L^2(\Omega)}+(\mathfrak{u}z,v)_{L^2(\Omega)}=(\mathfrak{f},v)_{L^2(\Omega)} \quad \forall v\in H_0^1(\Omega).
\end{equation}
Lax--Milgram Theorem immediately yields the well-posedness of problem \eqref{eq:weak_eq}. In particular, we have the following stability estimate
$
\|\nabla z\|_{L^2(\Omega)}\lesssim \|\mathfrak{f}\|_{L^2(\Omega)}.
$

%%%%%%%%%%%%%%%%%%%%%%%%%%%%%%%%%%%%%%%%%%%%%%%%%%%%%%%%%%%%
%%%%%%%%%%%%%%%%%%%%%%%%%%%%%%%%%%%%%%%%%%%%%%%%%%%%%%%%%%%%
%%%%%%%%%%%%%%%%%%%%%%%%%%%%%%%%%%%%%%%%%%%%%%%%%%%%%%%%%%%%
%%%%%%%%%%%%%%%%%%%%%%%%%%%%%%%%%%%%%%%%%%%%%%%%%%%%%%%%%%%%

\subsection{Finite Element Approximation}

In this section, we introduce a basic finite element approximation for the weak problem \eqref{eq:weak_eq} and review basic a posteriori error estimates. To accomplish this task, we first introduce some terminology and \FF{further} basic ingredients.

We denote by $\mathscr{T}=\{T\}$ a conforming partition of $\overline{\Omega}$ into simplices $T$ with size $h_T=\text{diam}(T)$ and define $h_{\mathscr{T}}:=\max_{T\in\mathscr{T}}h_T$. We denote by $\mathscr{S}$ the set of \emph{internal} $(d-1)$-dimensional interelement boundaries $S$ of $\mathscr{T}$. For $T \in \mathscr{T}$, we let $\mathscr{S}_T$ denote the subset of $\mathscr{S}$ which contains the sides \FF{of the element} $T$. We denote by $\FF{\mathcal{N}_S \subset \mathscr{T}}$ the subset that contains the two elements that have $S$ as a side, \FF{namely}, $\mathcal{N}_S=\{T^+,T^-\}$, where $T^+, T^- \in \mathscr{T}$ are such that $S = T^+ \cap T^-$. For $T \in \mathscr{T}$, we define the \emph{star} associated with the element $T$ as
\begin{equation}\label{def:patch}
\mathcal{N}_T:= \left \{ T^{\prime}\in\mathscr{T}: \mathscr{S}_{T}\cap \mathscr{S}_{T^\prime}\neq\emptyset \right \}.
\end{equation}
In an abuse of notation, below we denote by $\mathcal{N}_T$ either the set itself or the union of its elements. 

\EO{We define, for $T \in \mathscr{T}$, the \emph{shape coefficient} $\sigma_T$ of $T$ as the ratio of the diameter, i.e., $\mathrm{diam}(T) = h_T$, and the inball diameter of $T$, i.e., $2 \sup \{ r> 0: B_r(x) \subset T \textrm{ for } x \in T \}$. The \emph{shape coefficient} of a triangulation $\mathscr{T}$ corresponds to the quantity $\sigma_{\mathscr{T}}:= \max \{ T : T \in \mathscr{T}\}$. A sequence of triangulations $\mathbb{T} = \{ \mathscr{T} \}$, that is obtained by subsequent refinements of an initial mesh $\mathscr{T}_0$, is \emph{shape regular} if $\sup \{ \sigma_{\mathscr{T}}: \mathscr{T} \in \mathbb{T} \} \leq C$; see \FF{\cite[section 3.2.1]{Nochetto_etal2009}} for details.}
% ...
% , we denote by $\mathbb{T}$ the collection of conforming and shape regular meshes \FF{\cite[section 3.2.1]{Nochetto_etal2009}} that are obtained by subsequent refinements of an initial mesh $\mathscr{T}_0$.}

Given a mesh $\mathscr{T} \in \mathbb{T}$, we define the finite element space of continuous piecewise \FF{linear functions} as
\begin{equation}\label{def:piecewise_linear_set}
\mathbb{V}(\T):=\{v_\T\in C(\overline{\Omega}): v_\T|_T \in \mathbb{P}_{1}(T) \ \forall  T\in \T\}\cap H_0^1(\Omega).
\end{equation}

Given a discrete function $v_{\T} \in \mathbb{V}(\T)$, we define, for any internal side $S \in \Sides$, the jump or interelement residual  $\llbracket \nabla v_\mathscr{T}\cdot \boldsymbol{\nu} \rrbracket$ by
\[
\llbracket \nabla v_\mathscr{T}\cdot \boldsymbol{\nu} \rrbracket:= \boldsymbol{\nu}^{+} \cdot \nabla v_{\mathscr{T}}|^{}_{T^{+}} + \boldsymbol{\nu}^{-} \cdot \nabla v_{\mathscr{T}}|^{}_{T^{-}},
\]
where $\boldsymbol{\nu}^{+}, \boldsymbol{\nu}^{-}$ denote the unit normals to $S$ pointing towards $T^{+}$, $T^{-} \in \T$, respectively. Here, $T^{+}$, $T^{-} \in \T$ are such that $T^{+} \neq T^{-}$ and $\partial T^{+} \cap \partial T^{-} = S$.

With these ingredients at hand, we introduce a Galerkin approximation to problem \eqref{eq:weak_eq} as follows:
\begin{equation}\label{eq:discrete_eq}
z_\T\in\mathbb{V}(\T):
\quad
(\nabla z_\T, \nabla v_{\T})_{L^2(\Omega)} + (\mathfrak{u}z_\T,v_{\T})_{L^2(\Omega)} = (\mathfrak{f},v_{\T})_{L^2(\Omega)}
\end{equation}
for all $v^{}_\T\in \mathbb{V}(\T)$. Here, $\mathfrak{f}\in L^2(\Omega)$ and $\mathfrak{u}\in\mathbb{U}_{ad}$. The existence and uniqueness of a solution $z_{\T}\in  \mathbb{V}(\T)$ of problem \eqref{eq:discrete_eq} is standard. In particular, we have the stability estimate $\|\nabla z_\T\|_{L^2(\Omega)}\lesssim \|\mathfrak{f}\|_{L^2(\Omega)}$.

%%%%%%%%%%%%%%%%%%%%%%%%%%%%%%%%%%%%%%%%%%%%%%%%%%%%%%%%%%%%
%%%%%%%%%%%%%%%%%%%%%%%%%%%%%%%%%%%%%%%%%%%%%%%%%%%%%%%%%%%%
%%%%%%%%%%%%%%%%%%%%%%%%%%%%%%%%%%%%%%%%%%%%%%%%%%%%%%%%%%%%
%%%%%%%%%%%%%%%%%%%%%%%%%%%%%%%%%%%%%%%%%%%%%%%%%%%%%%%%%%%%

\subsection{An a Posteriori Error Estimate for the State Equation}

We introduce the following local error indicators and a posteriori error estimator associated to the discretization \eqref{eq:discrete_eq} of problem \eqref{eq:weak_eq}:
\begin{equation*}
\mathcal{E}_{T}^2:=h_{T}^2\|\mathfrak{f}-\mathfrak{u}z_{\T}\|_{L^2(T)}^2+h_{T}\|\llbracket \nabla z_{\T}\cdot \boldsymbol\nu \rrbracket\|_{L^2(\partial T\setminus\partial\Omega)}^2,\qquad \mathcal{E}_{\T}^2:= \sum_{T\in\T}\mathcal{E}_{T}^2.
\end{equation*}

We present the following global reliability result.

\begin{theorem}[global reliability of $\mathcal{E}$]\label{thm:global_reli_weak}
Let $\mathfrak{f} \in L^2(\Omega)$ and $\mathfrak{u}\in\mathbb{U}_{ad}$ be given. Let $z\in H_0^1(\Omega)$ be the unique solution to problem \eqref{eq:weak_eq} and let $z_\T\in\mathbb{V}(\T)$ be its finite element approximation obtained as the solution to \eqref{eq:discrete_eq}. We thus have
\[
\| \nabla(z-z_\T) \|_{L^2(\Omega)} \lesssim \mathcal{E}_{\T},
\]
with a hidden constant that is independent of $z$, $z_{\T}$, the size of the elements in $\T$, and $\#\T$ \EO{but depends on the shape coefficient of the triangulation $\T$, i.e., $\sigma_{\T}$, and the dimension $d$.}
\end{theorem}
\proof
Since $z$ solves \eqref{eq:weak_eq}, we invoke Galerkin orthogonality and an elementwise integration by parts formula to arrive at
\begin{multline*}
(\nabla (z-z_\T),\nabla v)_{L^2(\Omega)}+
(\mathfrak{u}(z-z_\T),v)_{L^2(\Omega)}
\\
= \sum_{T\in\T}
\int_{T}(\mathfrak{f}-\mathfrak{u}z_{\T})(v-I_\T v ) \mathrm{d} x
+\sum_{S\in\Sides}\int_{S} \llbracket\nabla z_\T\cdot\boldsymbol{\nu}\rrbracket ( v-I_\T v) \mathrm{d}x.
\end{multline*}
Here, $v \in H_0^1(\Omega)$ and $I_\T:L^1(\Omega)\rightarrow \mathbb{V}(\T)$ denotes the Cl\'ement interpolation operator \cite{MR2373954,MR0520174}. Standard approximation properties for $I_\T$ and the finite overlapping property of stars allow us to derive
\begin{multline*}
(\nabla (z-z_\T),\nabla v)_{L^2(\Omega)}+
(\mathfrak{u}(z-z_\T),v)_{L^2(\Omega)}
\\
\lesssim
\left[\sum_{T\in\T}h_T^2\|\mathfrak{f}-\mathfrak{u}z_\T\|_{L^2(T)}^2+h_T\|\llbracket\nabla z_\T\cdot\boldsymbol{\nu}\rrbracket\|_{L^2(\partial T\setminus \partial\Omega)}^2\right]^{\tfrac{1}{2}}\|\nabla v\|_{L^2(\Omega)}.
\end{multline*}
Set $v=z-z_\T\in H_0^1(\Omega)$ and use the fact that $\mathfrak{u}>0$ to conclude.
\qed

%%%%%%%%%%%%%%%%%%%%%%%%%%%%%%%%%%%%%%%%%%%%%%%%%%%%%%%%%%%%
%%%%%%%%%%%%%%%%%%%%%%%%%%%%%%%%%%%%%%%%%%%%%%%%%%%%%%%%%%%%
%%%%%%%%%%%%%%%%%%%%%%%%%%%%%%%%%%%%%%%%%%%%%%%%%%%%%%%%%%%%
%%%%%%%%%%%%%%%%%%%%%%%%%%%%%%%%%%%%%%%%%%%%%%%%%%%%%%%%%%%%
%%%%%%%%%%%%%%%%%%%%%%%%%%%%%%%%%%%%%%%%%%%%%%%%%%%%%%%%%%%%
%%%%%%%%%%%%%%%%%%%%%%%%%%%%%%%%%%%%%%%%%%%%%%%%%%%%%%%%%%%%
%%%%%%%%%%%%%%%%%%%%%%%%%%%%%%%%%%%%%%%%%%%%%%%%%%%%%%%%%%%%
%%%%%%%%%%%%%%%%%%%%%%%%%%%%%%%%%%%%%%%%%%%%%%%%%%%%%%%%%%%%

\section{The Optimal Control Problem}\label{sec:the_ocp}

In this section, we follow \cite[section 2]{MR2536007} and introduce a weak formulation for the optimal control problem \eqref{eq:minimize_cost_func}--\eqref{def:box_constraints}. In addition, we review first and second order optimality conditions and  introduce finite element discretization schemes.

%%%%%%%%%%%%%%%%%%%%%%%%%%%%%%%%%%%%%%%%%%%%%%%%%%%%%%%%%%%%
%%%%%%%%%%%%%%%%%%%%%%%%%%%%%%%%%%%%%%%%%%%%%%%%%%%%%%%%%%%%
%%%%%%%%%%%%%%%%%%%%%%%%%%%%%%%%%%%%%%%%%%%%%%%%%%%%%%%%%%%%
%%%%%%%%%%%%%%%%%%%%%%%%%%%%%%%%%%%%%%%%%%%%%%%%%%%%%%%%%%%%

\subsection{Weak Formulation and Existence of a Solution}\label{sec:weak_form}

Let $J$ be the cost functional defined in \eqref{def:cost_functional}. We consider the following weak version of the optimization problem \eqref{eq:minimize_cost_func}--\eqref{def:box_constraints}: Find
\begin{equation}\label{eq:minimize_weak}
\min \{ J(y,u): (y,u) \in H_0^1(\Omega)\times\mathbb{U}_{ad} \}
\end{equation}
subject to the state equation
\begin{equation}\label{eq:weak_state_eq}
(\nabla y,\nabla v)_{L^2(\Omega)}+(uy,v)_{L^2(\Omega)}=(f,v)_{L^2(\Omega)} \quad \forall v\in H_0^1(\Omega).
\end{equation}

The existence of an optimal solution $(\bar{y},\bar{u})\in H_0^1(\Omega)\times \mathbb{U}_{ad}$ for problem \eqref{eq:minimize_weak}--\eqref{eq:weak_state_eq} follows standard arguments; see \cite[Proposition 2.3]{MR2536007}.

%%%%%%%%%%%%%%%%%%%%%%%%%%%%%%%%%%%%%%%%%%%%%%%%%%%%%%%%%%%%
%%%%%%%%%%%%%%%%%%%%%%%%%%%%%%%%%%%%%%%%%%%%%%%%%%%%%%%%%%%%
%%%%%%%%%%%%%%%%%%%%%%%%%%%%%%%%%%%%%%%%%%%%%%%%%%%%%%%%%%%%
%%%%%%%%%%%%%%%%%%%%%%%%%%%%%%%%%%%%%%%%%%%%%%%%%%%%%%%%%%%%

\subsection{Optimality Conditions}\label{sec:opt_cond}

\FF{Due to the fact that} the optimal control problem \eqref{eq:minimize_weak}--\eqref{eq:weak_state_eq} is not convex, we discuss optimality conditions \FF{under the framework} of local solutions in $L^2(\Omega)$. \FF{To be precise,} a control $\bar{u} \in \mathbb{U}_{ad}$ is said to be locally optimal in $L^2(\Omega)$ for \eqref{eq:minimize_weak}--\eqref{eq:weak_state_eq} if there exists \FF{a constant} $\delta>0$ such that
$
J(\bar{y},\bar{u})\leq J(y,u)
$
for all $u\in\mathbb{U}_{ad}$ such that $\|u-\bar{u}\|_{L^2(\Omega)}\leq \delta$. Here, $\bar{y}$ and $y$ denote the states associated to $\bar{u}$ and $u$, respectively.

Let us introduce the set $\mathcal{U}:=\{u\in L^\infty(\Omega): \exists c > 0 \text{ such that } u(x) > c > 0 \text{ a.e. } x \in \Omega\}$. We immediately notice that $\mathbb{U}_{ad}\subset \mathcal{U}$. Having defined $\mathcal{U}$, we introduce the control-to-state map $\mathcal{S}$ \FF{as follows:} given a control $u\in \mathcal{U}$, \FF{$\mathcal{S}$} associates to it a unique state $y=\mathcal{S}u \in H_0^1(\Omega)$ \FF{solving} \eqref{eq:weak_state_eq}. With these ingredients at hand, we define the reduced cost functional $j: \mathcal{U}\to \mathbb{R}_{0}^{+}$ by
\begin{equation*}
j(u)=J(\mathcal{S}u,u):=\frac{1}{2}\|\mathcal{S}u-y^{}_\Omega\|_{L^2(\Omega)}^2+\frac{\FF{\alpha}}{2}\|u\|_{L^2(\Omega)}^2.
\end{equation*}

\FF{We are now in position to} formulate first order optimality conditions: if $\bar{u}$ is locally optimal for problem \eqref{eq:minimize_weak}--\eqref{eq:weak_state_eq}, then \cite[Proposition 2.10]{MR2536007}
\begin{equation}\label{eq:gat_deri}
j'(\bar{u})(u-\bar{u})\geq 0 \quad \forall u\in\mathbb{U}_{ad}.
\end{equation}
Here, $j'(\bar{u})$ denotes the Gate\^aux derivative of the functional $j$ at $\bar{u}$ in the direction $u-\bar{u}$. \EO{We notice that, for $u,v\in \mathbb{U}_{ad}$, $j'(u)v = (\alpha u - yp,v)_{L^2(\Omega)}$ \cite[equation (2.6)]{MR2536007} and immediately comment that $\mathcal{S}$ and $j$ are not Fr\'echet differentiable with respect to the $L^2(\Omega)$-topology \cite[Remark 2.8]{MR2536007}.} To \FF{investigate} the inequality \eqref{eq:gat_deri}, we introduce the \emph{adjoint variable} $p\in H_0^1(\Omega)$ as the unique solution to the \emph{adjoint equation}
\begin{equation}\label{eq:adj_eq}
(\nabla w,\nabla p)_{L^2(\Omega)}+(up,w)_{L^2(\Omega)}=(y-y_\Omega,w)_{L^2(\Omega)} \quad \forall  w\in H_0^1(\Omega),
\end{equation}
where $y = \mathcal{S}u$ solves \eqref{eq:weak_state_eq}. Observe that problem \eqref{eq:adj_eq} is well-posed.

\FF{With the previous ingredients at hand, we} reformulate first order optimality conditions as follows; see \cite[Proposition 2.10 and equation (2.6)]{MR2536007}.
\begin{theorem}[first order optimality conditions]
Every locally optimal control $\bar{u}\in\mathbb{U}_{ad}$ for problem
\eqref{eq:minimize_weak}--\eqref{eq:weak_state_eq} satisfies, together with the state $\bar{y}\in H_0^1(\Omega)$ and the adjoint state $\bar{p}\in H_0^1(\Omega)$, the variational inequality
\begin{equation}\label{eq:var_ineq}
(\alpha\bar{u}-\bar{y}\bar{p},u-\bar{u})_{L^2(\Omega)}\geq 0 \quad \forall u\in\mathbb{U}_{ad}.
\end{equation}
Here, $\bar{p}$ denotes the solution to \eqref{eq:adj_eq} with $y$ replaced by $\bar{y}=\mathcal{S}\bar{u}$.
\label{thm:1st_order_cond}
\end{theorem}

Let us now introduce the projection operator $\Pi_{[\texttt{a},\texttt{b}]} : L^1(\Omega) \rightarrow  \mathbb{U}_{ad}$ as 
\begin{equation}\label{def:projector_pi}
\Pi_{[\texttt{a},\texttt{b}]}(v) := \min\{ \texttt{b}, \max\{ v, \texttt{a}\} \} \textrm{ a.e. in } \Omega.
\end{equation} 
\FF{This operator allows us to} present the following projection formula \cite[equation (2.7)]{MR2536007}: If $\bar u$ denotes a locally optimal control for \eqref{eq:minimize_weak}--\eqref{eq:weak_state_eq}, then
\begin{equation}\label{eq:projection_control}
\bar{u}(x):=\Pi_{[\texttt{a},\texttt{b}]}(\alpha^{-1}\bar{y}(x)\bar{p}(x)) \textrm{ a.e.}~x \in \Omega.
\end{equation}

Let $\bar{u}\in\mathbb{U}_{ad}$ be a control that satisfies the first order necessary optimality condition \eqref{eq:var_ineq}. In what follows, we will assume that there exists a constant $\mu > 0$ such that
\begin{equation}\label{eq:second_order_cond}
j''(\bar{u})v^2 \geq \mu \|v\|_{L^2(\Omega)}^2 \quad \forall v\in L^\infty(\Omega);
\end{equation}
see \cite[Assumption 2.20]{MR2536007}. \FF{Here, for each $v\in L^{\infty}(\Omega)$, we have that
\begin{equation*}
j''(\bar{u})v^2 := \|\mathcal{S}'(\bar{u})v\|_{L^2(\Omega)}^2 + (\bar{y} - y_{\Omega}, \mathcal{S}''(\bar{u})v^2)_{L^2(\Omega)} + \alpha \|v\|_{L^2(\Omega)}^2,
\end{equation*}
where $\mathcal{S}'(\bar{u})v$ and $\mathcal{S}''(\bar{u})v^2$ are defined as in \cite[Lemma 2.9]{MR2536007}.} We notice that assumption \eqref{eq:second_order_cond} is fulfilled if $\|\mathcal{S}\bar{u}-y_\Omega\|_{L^2(\Omega)}$ is sufficiently small or $\alpha>0$ is sufficiently large; see \cite[Remark 2.21]{MR2536007} for details.

The following result states that every control $\bar{u}$ that satisfies 
\eqref{eq:var_ineq} and \eqref{eq:second_order_cond} is a local solution for problem \eqref{eq:minimize_weak}--\eqref{eq:weak_state_eq}; see \cite[Theorem 2.24]{MR2536007}.

\begin{theorem}[local optimality]\label{thm:local_opt}
Let $\bar{u}\in\mathbb{U}_{ad}$ be a local solution to \eqref{eq:minimize_weak}--\eqref{eq:weak_state_eq} satisfying 
the necessary and sufficient optimality conditions \eqref{eq:var_ineq} and \eqref{eq:second_order_cond}. Then, there exist positive constants $\delta,\sigma>0$ such that
\begin{equation*}
j(u)\geq j(\bar{u})+\sigma\|u-\bar{u}\|_{L^2(\Omega)}^2 \quad \forall u\in\mathbb{U}_{ad}\cap B_{\delta}(\bar{u}),
\end{equation*}
where $B_{\delta}(\bar{u})$ denotes the closed ball in $L^2(\Omega)$ with center at $\bar{u}$ and radius $\delta$.
\end{theorem}

We conclude this section with the following estimate \cite[Proposition 2.22]{MR2536007}: Let $u ,v \in \mathbb{U}_{ad}$ and $w \in L^\infty(\Omega)$. Then, there exists $\mathfrak{C} > 0$, depending on $\|f\|_{L^2(\Omega)}$ and $\|y_\Omega\|_{L^2(\Omega)}$, such that 
\begin{equation}\label{eq:lipschitz_j''}
|j''(u)w^2-j''(v)w^2|\leq \mathfrak{C}\|u - v \|_{L^2(\Omega)}\|w\|_{L^2(\Omega)}^2.
\end{equation}

%%%%%%%%%%%%%%%%%%%%%%%%%%%%%%%%%%%%%%%%%%%%%%%%%%%%%%%%%%%%
%%%%%%%%%%%%%%%%%%%%%%%%%%%%%%%%%%%%%%%%%%%%%%%%%%%%%%%%%%%%
%%%%%%%%%%%%%%%%%%%%%%%%%%%%%%%%%%%%%%%%%%%%%%%%%%%%%%%%%%%%
%%%%%%%%%%%%%%%%%%%%%%%%%%%%%%%%%%%%%%%%%%%%%%%%%%%%%%%%%%%%

\subsection{Finite Element Approximation}

In this section, we introduce two finite element discretization schemes for our optimal control problem. 

%%%%%%%%%%%%%%%%%%%%%%%%%%%%%%%%%%%%%%%%%%%%%%%%%%%%%%%%%%%%
%%%%%%%%%%%%%%%%%%%%%%%%%%%%%%%%%%%%%%%%%%%%%%%%%%%%%%%%%%%%
%%%%%%%%%%%%%%%%%%%%%%%%%%%%%%%%%%%%%%%%%%%%%%%%%%%%%%%%%%%%
%%%%%%%%%%%%%%%%%%%%%%%%%%%%%%%%%%%%%%%%%%%%%%%%%%%%%%%%%%%%

\subsubsection{The Fully Discrete Scheme}\label{sec:fully_discrete_scheme}

To approximate a control variable, we introduce the space of piecewise constant functions
\begin{equation*}
\mathbb{U}(\mathscr{T}):=\{ u_\T\in L^\infty(\Omega): u_\T|_T\in \mathbb{P}_0(T) \ \forall  T\in \T\}
\end{equation*}
and define the discrete admissible set $\mathbb{U}_{ad}(\mathscr{T}):=\mathbb{U}(\mathscr{T})\cap \mathbb{U}_{ad}$. \EO{The} state and adjoint state \FF{variables, associated to a locally optimal control,} are discretized by using the finite element space $\mathbb{V}(\T)$ defined in \eqref{def:piecewise_linear_set}. With this setting at hand, the fully discrete scheme reads as follows: Find $\min J(y_{\T},u_{\T})$ subject to the discrete state equation
\begin{equation}
\label{eq:discrete_state_equation}
y^{}_\mathscr{T} \in \mathbb{V}(\mathscr{T}): \quad
(\nabla y^{}_\mathscr{T},\nabla v^{}_\mathscr{T})_{L^2(\Omega)}+(u_{\T}y_{\T},v_\T)_{L^2(\Omega)}  =  (f,v^{}_\mathscr{T})^{}_{L^2(\Omega)}
\end{equation}
for all $v^{}_\mathscr{T} \in \mathbb{V}(\mathscr{T})$ and the discrete constraints $u^{}_{\mathscr{T}} \in \mathbb{U}_{ad}(\T)$. The fully discrete scheme admits at least a solution; see \cite[Section 3]{MR2536007} for details. In addition, if $\bar u^{}_{\T}$ denotes a discrete local solution, then
\begin{equation*}
\label{eq:discrete_opt_system}
 (\alpha\bar{u}^{}_\mathscr{T}-\bar{y}_{\T}\bar{p}_{\T},u^{}_\mathscr{T}-\bar{u}^{}_\mathscr{T})^{}_{L^2(\Omega)}  \geq  0 \quad \forall u^{}_\mathscr{T} \in \mathbb{U}_{ad}(\T),
\end{equation*}
where $\bar{p}^{}_\mathscr{T} \in \mathbb{V}(\T)$ is such that
\begin{equation}
\label{eq:discrete_adjoint_equation}
(\nabla w^{}_\mathscr{T},\nabla \bar{p}^{}_\mathscr{T})_{L^2(\Omega)}+(\bar{u}_{\T}\bar{p}_\T,w_\T)_{L^2(\Omega)}    =   (\bar{y}^{}_\mathscr{T}-y^{}_{\Omega},w^{}_\mathscr{T})^{}_{L^2(\Omega)} 
\end{equation}
for all $w^{}_\mathscr{T} \in \mathbb{V}(\T)$; see \cite[equations (3.6) and (3.7)]{MR2536007}.

%%%%%%%%%%%%%%%%%%%%%%%%%%%%%%%%%%%%%%%%%%%%%%%%%%%%%%%%%%%%
%%%%%%%%%%%%%%%%%%%%%%%%%%%%%%%%%%%%%%%%%%%%%%%%%%%%%%%%%%%%
%%%%%%%%%%%%%%%%%%%%%%%%%%%%%%%%%%%%%%%%%%%%%%%%%%%%%%%%%%%%
%%%%%%%%%%%%%%%%%%%%%%%%%%%%%%%%%%%%%%%%%%%%%%%%%%%%%%%%%%%%

\subsubsection{The Semi-discrete Scheme}\label{sec:semi_discrete_scheme}
In this section, we introduce the so-called variational discretization approach for \eqref{eq:minimize_weak}--\eqref{eq:weak_state_eq}. This scheme discretizes only the state space; the control space $\mathbb{U}_{ad}$ is not discretized. The scheme induces a discretization of an optimal control variable by projecting, \FF{in view of the operator introduced in \eqref{def:projector_pi},} an optimal discrete adjoint state into $\mathbb{U}_{ad}$. The semi-discrete scheme is defined as follows: Find $\min J(y_{\T},\mathsf{u})$ subject to the discrete state equation
\begin{equation}
\label{eq:semi_discrete_state_equation}
y^{}_\mathscr{T} \in \mathbb{V}(\mathscr{T}): \quad
(\nabla y^{}_\mathscr{T},\nabla v^{}_\mathscr{T})_{L^2(\Omega)}+(\mathsf{u}y_{\T},v_\T)_{L^2(\Omega)}  =  (f,v^{}_\mathscr{T})^{}_{L^2(\Omega)}
\end{equation}
for all $v^{}_\mathscr{T} \in \mathbb{V}(\mathscr{T})$ and the constraints $\mathsf{u} \in \mathbb{U}_{ad}$. As in the fully discrete case, this problem admits at least a solution and, if $\bar{\mathsf{u}}$ denotes a local solution, then
\begin{equation*}
\label{eq:semi_var_ineq}
 (\alpha\bar{\mathsf{u}}-\bar{y}_{\T}\bar{p}_{\T},u-\bar{\mathsf{u}})^{}_{L^2(\Omega)}  \geq  0 \quad \forall u \in \mathbb{U}_{ad},
\end{equation*}
where $\bar{p}^{}_\mathscr{T} \in \mathbb{V}(\T)$ solves 
\begin{equation}
\label{eq:semi_discrete_adjoint_equation}
(\nabla w^{}_\mathscr{T},\nabla \bar{p}^{}_\mathscr{T})_{L^2(\Omega)}+(\bar{\mathsf{u}}\bar{p}_\T,w_\T)_{L^2(\Omega)}    =   (\bar{y}^{}_\mathscr{T}-y^{}_{\Omega},w^{}_\mathscr{T})^{}_{L^2(\Omega)} 
\end{equation}
for all $w^{}_\mathscr{T} \in \mathbb{V}(\T)$. Here, $\bar{y}^{}_\mathscr{T}=\bar{y}^{}_\mathscr{T}(\bar{\mathsf{u}})$ solves  \eqref{eq:semi_discrete_state_equation} with $\mathsf{u}=\bar{\mathsf{u}}$.

%%%%%%%%%%%%%%%%%%%%%%%%%%%%%%%%%%%%%%%%%%%%%%%%%%%%%%%%%%%%
%%%%%%%%%%%%%%%%%%%%%%%%%%%%%%%%%%%%%%%%%%%%%%%%%%%%%%%%%%%%
%%%%%%%%%%%%%%%%%%%%%%%%%%%%%%%%%%%%%%%%%%%%%%%%%%%%%%%%%%%%
%%%%%%%%%%%%%%%%%%%%%%%%%%%%%%%%%%%%%%%%%%%%%%%%%%%%%%%%%%%%
%%%%%%%%%%%%%%%%%%%%%%%%%%%%%%%%%%%%%%%%%%%%%%%%%%%%%%%%%%%%
%%%%%%%%%%%%%%%%%%%%%%%%%%%%%%%%%%%%%%%%%%%%%%%%%%%%%%%%%%%%
%%%%%%%%%%%%%%%%%%%%%%%%%%%%%%%%%%%%%%%%%%%%%%%%%%%%%%%%%%%%
%%%%%%%%%%%%%%%%%%%%%%%%%%%%%%%%%%%%%%%%%%%%%%%%%%%%%%%%%%%%

\section{A Posteriori Error Analysis for the Fully Discrete Scheme}\label{sec:a_posteriori_fully}

In this section, we devise and analyze an a posteriori error estimator for the fully discrete scheme. The error estimator will be \FF{formed by} the sum of three contributions: two contributions related to the discretization of the state and adjoint equations and a \FF{one contribution} associated to the discretization of the admissible control set $\mathbb{U}_{ad}$.

To begin with our studies we introduce, on the basis of the projection operator $\Pi_{[\texttt{a},\texttt{b}]}$ defined in \eqref{def:projector_pi}, the auxiliary variable
\begin{equation}\label{def:control_tilde}
\tilde{u}:=\Pi_{[\texttt{a},\texttt{b}]}\left(\alpha^{-1}\bar{y}_{\T}\bar{p}_{\T}\right).
\end{equation}
A key property in favor of the definition of $\tilde u$ is that it satisfies the following variational inequality \cite[Lemma 2.26]{Troltzsch}:
\begin{equation}\label{eq:var_ineq_u_tilde}
(\alpha\tilde{u}-\bar{y}_{\T}\bar{p}_{\T},u-\tilde{u})_{L^2(\Omega)}\geq 0 \quad \forall u\in \mathbb{U}_{ad}.
\end{equation}

With the variable $\tilde{u}$ at hand, we present the following result which is instrumental for our a posteriori error analysis.

\begin{theorem}[auxiliary estimate]
\label{thm:error_bound_control_tilde}
Let $\bar u \in \mathbb{U}_{ad}$ be a local solution to \eqref{eq:minimize_weak}--\eqref{eq:weak_state_eq} satisfying the sufficient second order optimality condition \eqref{eq:second_order_cond}.
\FF{Let $\bar u_{\T}$ be a local minimum of the fully discrete optimal control problem with $\bar{y}_{\T}$ and $\bar{p}_\T$ being the corresponding state and adjoint state, respectively. If $\bar{y}_{\T}$ and $\bar{p}_\T$ satisfy, on the mesh $\mathscr{T}$, the bound}
\begin{equation}\label{eq:assumption_mesh}
\|\bar{y}\bar{p}-\bar{y}_\T\bar{p}_\T\|_{L^2(\Omega)} \leq \alpha\mu(2\mathfrak{C})^{-1},
\end{equation}
then
\begin{equation}\label{eq:inequality_control_tilde}
\frac{\mu}{2}\|\bar{u}-\tilde{u}\|_{L^2(\Omega)}^2
\leq
(j'(\tilde{u})-j'(\bar{u}))(\tilde{u}-\bar{u}).
\end{equation}
The constants $\mu$ and $\mathfrak{C}$ are given as in  \eqref{eq:second_order_cond} and \eqref{eq:lipschitz_j''}, respectively.
\end{theorem}
\proof
Since $\tilde{u}-\bar{u}\in L^\infty(\Omega)$, we are allowed to set $v=\tilde{u}-\bar{u}$ in the second order optimality condition \eqref{eq:second_order_cond}. This yields
\begin{equation}\label{eq:diff_1}
\mu\|\tilde{u}-\bar{u}\|_{L^2(\Omega)}^2 \leq j''(\bar{u})(\tilde{u}-\bar{u})^2.
\end{equation}
On the other hand, in view of the mean value theorem, we obtain
$
(j'(\tilde{u})-j'(\bar{u}))(\tilde{u}-\bar{u})=j''(\zeta)(\tilde{u}-\bar{u})^2,
$
where $\zeta=\bar{u}+\theta_{\T}(\tilde{u}-\bar{u})$ and $\theta_{\T} \in (0,1)$. Inequality \eqref{eq:diff_1} thus yields
\begin{align}\label{eq:ineq_u_tilde_bar}
\mu\|\tilde{u}-\bar{u}\|_{L^2(\Omega)}^2 
& \leq (j'(\tilde{u})-j'(\bar{u}))(\tilde{u}-\bar{u}) + (j''(\bar{u})-j''(\zeta))(\tilde{u}-\bar{u})^2.
\end{align}
Let us now concentrate \EO{on} the second term on the right-hand side of inequality \eqref{eq:ineq_u_tilde_bar}. To accomplish this task, let us first invoke \eqref{eq:lipschitz_j''} to arrive at
\begin{equation*}
(j''(\bar{u})-j''(\zeta))(\tilde{u}-\bar{u})^2 
\leq \mathfrak{C} \|\tilde{u}-\bar{u}\|_{L^{2}(\Omega)}\|\tilde{u}-\bar{u}\|_{L^2(\Omega)}^2,
\end{equation*}
where we have also used that $\theta_{\T} \in (0,1)$. Invoke \eqref{eq:projection_control} and \eqref{def:control_tilde}, the Lipschitz property of the projection operator $\Pi_{[\texttt{a},\texttt{b}]}$, and assumption \eqref{eq:assumption_mesh} to conclude
\begin{equation*}
(j''(\bar{u})-j''(\zeta))(\tilde{u}-\bar{u})^2 
\leq \mathfrak{C}\alpha^{-1} \|\bar{y}\bar{p} - \bar{y}_\T\bar{p}_\T \|_{L^{2}(\Omega)}\|\tilde{u}-\bar{u}\|_{L^2(\Omega)}^2
\leq \frac{\mu}{2}\|\tilde{u}-\bar{u}\|_{L^2(\Omega)}^2.
\end{equation*}
Replacing this inequality into \eqref{eq:ineq_u_tilde_bar} yields the desired inequality \eqref{eq:inequality_control_tilde}. 
\qed

\begin{remark}[a sufficient condition for estimate \eqref{eq:assumption_mesh}] An estimate that guarantees assumption \eqref{eq:assumption_mesh} reads as follows:
\[
\|\nabla(\bar{y} -\bar{y}_\T)\|_{L^2(\Omega)}+ \|\nabla(\bar{p} -\bar{p}_\T)\|_{L^2(\Omega)}\leq \alpha\mu(2\mathfrak{C}C_{cond})^{-1},
\]
where $C_{cond}$ is defined as in \eqref{eq:constant_suff}. In fact, let  $q \in \{2,4\}$ and let $C_{q}$ be the best constant associated to the embedding $H_0^1(\Omega)\hookrightarrow L^{q}(\Omega)$, i.e., $C_{q}$ is the best constant such that $\|v\|_{L^{q}(\Omega)}\leq C_{q}\|\nabla v\|_{L^2(\Omega)}$ for all $v \in H_0^1(\Omega)$. Hence, an application of the Cauchy--Schwarz inequality, the aforementioned Sobolev embeddings, and the stability of the adjoint and discrete state equations allow us to obtain the bounds
\begin{multline*}
\|\bar{y}\bar{p} - \bar{y}_\T\bar{p}_\T \|_{L^{2}(\Omega)}
\leq 
\|\bar{p}\|_{L^4(\Omega)}\|\bar{y} -\bar{y}_\T\|_{L^4(\Omega)}+\|\bar{y}_\T\|_{L^4(\Omega)}\|\bar{p} -\bar{p}_\T \|_{L^4(\Omega)}\\
\leq C_{4}^2
\left(\|\nabla\bar{p}\|_{L^2(\Omega)}\|\nabla(\bar{y} -\bar{y}_\T)\|_{L^2(\Omega)}+\|\nabla\bar{y}_\T\|_{L^2(\Omega)}\|\nabla(\bar{p} -\bar{p}_\T)\|_{L^2(\Omega)}\right)
\\
\leq  C_{4}^2
C_{2}
\left[( C_{2}^2 \|f\|_{L^2(\Omega)}+\|y_\Omega\|_{L^2(\Omega)})\|\nabla(\bar{y} -\bar{y}_\T)\|_{L^2(\Omega)}\right.\\
\left.+  \|f\|_{L^2(\Omega)}\|\nabla(\bar{p} -\bar{p}_\T)\|_{L^2(\Omega)}\right].
\end{multline*}
The desired result can thus be concluded by setting
\begin{equation}\label{eq:constant_suff}
C_{cond}:=C_{4}^2C_{2}\max\{C_{2}^2\|f\|_{L^2(\Omega)}+\|y_\Omega\|_{L^2(\Omega)},\|f\|_{L^2(\Omega)}\}.
\end{equation}
\end{remark}

%%%%%%%%%%%%%%%%%%%%%%%%%%%%%%%%%%%%%%%%%%%%%%%%%%%%%%%%%%%%
%%%%%%%%%%%%%%%%%%%%%%%%%%%%%%%%%%%%%%%%%%%%%%%%%%%%%%%%%%%%
%%%%%%%%%%%%%%%%%%%%%%%%%%%%%%%%%%%%%%%%%%%%%%%%%%%%%%%%%%%%
%%%%%%%%%%%%%%%%%%%%%%%%%%%%%%%%%%%%%%%%%%%%%%%%%%%%%%%%%%%%

\subsection{Global Reliability Analysis}
\label{sec:global_fully}

The goal of this section is to derive an upper bound for the corresponding total error in terms of a devised a posteriori error estimator. The analysis relies on estimates on the error between a solution to the fully discrete optimal control problem \eqref{eq:discrete_state_equation}--\eqref{eq:discrete_adjoint_equation} and auxiliary variables that we define in what follows. 

We first define the variable $\hat y \in H_{0}^{1}(\Omega)$ as the solution to
\begin{equation}\label{eq:hat_functions}
(\nabla\hat{y},\nabla v)_{L^2(\Omega)}+(\bar{u}_\T\hat{y},v)_{L^2(\Omega)} =  (f,v)^{}_{L^2(\Omega)} \quad \forall v \in H_0^1(\Omega).
\end{equation}

Define now, for $T\in\T$, the local error indicators 
\begin{equation}\label{eq:indicators_state_eq}
\mathcal{E}_{st,T}^2:=h_T^2\|f-\bar{u}_\T\bar{y}_\T\|_{L^2(T)}^2+h_T\|\llbracket\nabla \bar{y}_\T\cdot\boldsymbol{\nu}\rrbracket\|_{L^2(\partial T\setminus\partial\Omega)}^2,
\end{equation}
and the global a posteriori error estimator associated to the finite element discretization of the state equation 
\begin{equation}\label{eq:error_estimator_state_eq}
\mathcal{E}_{st,\T}:=\left(\sum_{T\in\T}\mathcal{E}_{st,T}^2\right)^{\frac{1}{2}}.
\end{equation}
An application of Theorem \ref{thm:global_reli_weak}, with $\mathfrak{f}=f$ and $\mathfrak{u}=\bar{u}_\T$, immediately yields the a posteriori error estimate
\begin{equation}\label{eq:estimate_state_hat_discrete_st}
\|\nabla (\hat{y}-\bar{y}_\T)\|_{L^2(\Omega)}
\lesssim \mathcal{E}_{st,\T}.
\end{equation}

Let $\hat p \in H_{0}^{1}(\Omega)$ be the solution to
\begin{equation}\label{eq:hat_functions_p}
(\nabla w,\nabla \hat{p})_{L^2(\Omega)}+\left(\bar{u}_\T\hat{p},w\right)_{L^2(\Omega)} =  (\bar{y}_\T-y^{}_{\Omega},w)^{}_{L^2(\Omega)} \quad \forall w\in H_0^1(\Omega).
\end{equation}
Define, for $T\in\T$, the local error indicators
\begin{equation}
\label{eq:error_indicator_adjoint_eq}
\mathcal{E}_{adj,T}^2:=h_T^2\|\bar{y}_\T-y_\Omega-\bar{u}_{\T}\bar{p}_{\T}\|_{L^2(T)}^2+h_T\|\llbracket\nabla \bar{p}_\T\cdot\boldsymbol{\nu}\rrbracket\|_{L^2(\partial T\setminus\partial\Omega)}^2,
\end{equation}
and the a posteriori error estimator associated to the discretization of the adjoint equation
\begin{equation}\label{eq:error_estimator_adjoint_eq}
\mathcal{E}_{adj,\T}:=\left(\sum_{T\in\T}\mathcal{E}_{adj,T}^2\right)^{\frac{1}{2}}.
\end{equation}
An application of Theorem \ref{thm:global_reli_weak}, again, with $\mathfrak{f}=\bar{y}_\T - y_\Omega$ and $\mathfrak{u}=\bar{u}_\T$, immediately yields the a posteriori error bound
\begin{equation}\label{eq:estimate_state_hat_discrete_adj}
\|\nabla (\hat{p}-\bar{p}_\T)\|_{L^2(\Omega)}
\lesssim \mathcal{E}_{adj,\T}.
\end{equation}

Finally, we introduce local error indicators and an a posteriori error estimator associated to the discretization of a control variable. To be precise, we define, on the basis  of the auxiliary variable $\tilde{u}$, defined in \eqref{def:control_tilde},
\begin{equation}\label{eq:error_estimator_control}
\mathcal{E}_{ct,T}^2:=\|\tilde{u}-\bar{u}_\T\|_{L^2(T)}^2,\qquad 
\mathcal{E}_{ct,\T}:=\left(\sum_{T\in\T}\mathcal{E}_{ct,T}^2\right)^\frac{1}{2}.
\end{equation}

After having defined error estimators associated to the discretization of the state and adjoint equations and the admissible control set, we introduce an a posteriori error estimator for the fully discrete problem that approximates solutions to problem \eqref{eq:minimize_weak}--\eqref{eq:weak_state_eq}. The error estimator can be decomposed as the sum of three contributions:
\begin{equation}\label{def:error_estimator_ocp}
\mathcal{E}_{ocp,\T}^2:=\mathcal{E}_{st,\T}^2 + \mathcal{E}_{adj,\T}^2 + \mathcal{E}_{ct,\T}^2.
\end{equation}
The estimators $\mathcal{E}_{st,\T}$, $\mathcal{E}_{adj,\T}$, and $\mathcal{E}_{ct,\T}$ are defined as in \eqref{eq:error_estimator_state_eq}, \eqref{eq:error_estimator_adjoint_eq}, and \eqref{eq:error_estimator_control}, respectively.

We are now ready to state and prove the main result of this section. As a final ingredient, we introduce $e_{\bar{y}} := \bar{y} - \bar{y}_{\T}$, $e^{}_{\bar{p}} := \bar{p} - \bar{p}_{\T}$, and $e_{\bar{u}} := \bar{u}- \bar{u}_{\T}$.

\begin{theorem}[global reliability]\label{thm:global_rel}
Let $\bar u \in \mathbb{U}_{ad}$ be a local solution to \eqref{eq:minimize_weak}--\eqref{eq:weak_state_eq} satisfying the sufficient second order condition \eqref{eq:second_order_cond}. Let $\bar u_{\T}$ be a local minimum of the fully discrete optimal control problem with $\bar{y}_{\T}$ and $\bar{p}_\T$ being the corresponding state and adjoint state, respectively. If \FF{$\bar{y}_{\T}$ and $\bar{p}_\T$ satisfy, on the mesh $\T$, the bound \eqref{eq:assumption_mesh}}, then
\begin{equation}\label{eq:global_rel}
\|\nabla{e}_{\bar{y}}\|^2_{L^2(\Omega)} + \|\nabla{e}_{\bar{p}}\|^2_{L^2(\Omega)} + \| e_{\bar{u}} \|^2_{L^2(\Omega)}
\lesssim 
\mathcal{E}_{ocp,\T}^2,
\end{equation}
with a hidden constant that is independent of continuous and discrete optimal variables, the size of the elements in $\T$, and $\#\T$.
\end{theorem}
\proof
We proceed on the basis of four steps.

\underline{Step 1.} The goal of this step is to control the term $\|e_{\bar{u}}\|_{L^2(\Omega)}$. To accomplish this task, we invoke the auxiliary variable $\tilde{u}=\Pi_{[\texttt{a},\texttt{b}]}(\alpha^{-1}\bar{y}_{\T}\bar{p}_{\T})$, a triangle inequality, and both definitions in \eqref{eq:error_estimator_control} to arrive at
\begin{equation}\label{eq:step1_i}
\|e_{\bar{u}}\|_{L^2(\Omega)}\leq \|\bar{u}-\tilde{u}\|_{L^2(\Omega)}+\mathcal{E}_{ct,\T}.
\end{equation}
Let us now concentrate \EO{on} the first term on the right-hand side of the previous estimate. Set $u=\tilde{u}$ in \eqref{eq:gat_deri} and $u=\bar{u}$ in \eqref{eq:var_ineq_u_tilde} to obtain, respectively,
\begin{equation*}
j'(\bar{u})(\tilde{u}-\bar{u})\geq 0,
\qquad
(\alpha\tilde{u}-\bar{p}_\T\bar{y}_\T,\bar{u}-\tilde{u})_{L^2(\Omega)}\geq 0.
\end{equation*}
With these estimates at hand, we invoke inequality \eqref{eq:inequality_control_tilde} to arrive at
\begin{align}\label{eq:step_ineq_1}
\frac{\mu}{2}\|\bar{u}-\tilde{u}\|_{L^2(\Omega)}^2
&\leq 
j'(\tilde{u})(\tilde{u}-\bar{u})-j'(\bar{u})(\tilde{u}-\bar{u})
\leq
j'(\tilde{u})(\tilde{u}-\bar{u})\\ \nonumber
&=
(\alpha\tilde{u}-\tilde{p}\tilde{y},\tilde{u}-\bar{u})_{L^2(\Omega)}
\leq
(\bar{p}_\T\bar{y}_\T-\tilde{p}\tilde{y},\tilde{u}-\bar{u})_{L^2(\Omega)},
\end{align}
where the auxiliary variables $\tilde{y}, \tilde{p}\in H_0^1(\Omega)$ are defined as follows:
\begin{equation}\label{eq:tilde_y}
\tilde{y}:
\quad
(\nabla\tilde{y},\nabla v)_{L^2(\Omega)}+(\tilde{u}\tilde{y},v)_{L^2(\Omega)} = (f,v)^{}_{L^2(\Omega)} \quad \forall v \in H_0^1(\Omega),
\end{equation}
and 
\begin{equation*}
\tilde{p}:
\quad
(\nabla w,\nabla \tilde{p})_{L^2(\Omega)}+(\tilde{u}\tilde{p},w)_{L^2(\Omega)} =  (\tilde{y}-y^{}_{\Omega},w)^{}_{L^2(\Omega)} \quad 
\forall w\in H_0^1(\Omega),
\end{equation*}
respectively. 

Adding and subtracting the term $\tilde{p}\bar{y}_\T$ in inequality \eqref{eq:step_ineq_1} and utilizing a generalized H\"older's inequality we arrive at
\begin{multline*}
\|\bar{u}-\tilde{u}\|_{L^2(\Omega)}^2
\lesssim 
([\bar{p}_\T-\tilde{p}]\bar{y}_\T+\tilde{p}[\bar{y}_\T-\tilde{y}],\tilde{u}-\bar{u})_{L^2(\Omega)}\\
\lesssim (\|\bar{p}_\T-\tilde{p}\|_{L^4(\Omega)}\|\bar{y}_\T\|_{L^4(\Omega)}+\|\bar{y}_\T-\tilde{y}\|_{L^4(\Omega)}\|\tilde{p}\|_{L^4(\Omega)})\|\tilde{u}-\bar{u}\|_{L^2(\Omega)}.
\end{multline*}
In view of the Sobolev embedding $H_0^1(\Omega)\hookrightarrow L^4(\Omega)$, we thus conclude that
\begin{multline}\label{eq:step1_ii}
\|\bar{u}-\tilde{u}\|_{L^2(\Omega)}
\lesssim 
\|\nabla(\bar{p}_\T-\tilde{p})\|_{L^2(\Omega)}\|\nabla\bar{y}_\T\|_{L^2(\Omega)}\\ +\|\nabla(\bar{y}_\T-\tilde{y})\|_{L^2(\Omega)}\|\nabla\tilde{p}\|_{L^2(\Omega)}.
\end{multline}
We now control $\|\nabla\bar{y}_\T\|_{L^2(\Omega)}$ and $\|\nabla\tilde{p}\|_{L^2(\Omega)}$. Set $v_\T=\bar{y}_\T$ in \eqref{eq:discrete_state_equation} and use that $\bar{u}_\T\in\mathbb{U}_{ad}$ to obtain $\|\nabla\bar{y}_\T\|_{L^2(\Omega)}\lesssim \|f\|_{L^2(\Omega)}$. Similar arguments yield
\begin{equation*}
\|\nabla\tilde{p}\|_{L^2(\Omega)}
\lesssim 
\|\tilde{y}-y_\Omega\|_{L^2(\Omega)}
\lesssim 
\|\nabla\tilde{y}\|_{L^2(\Omega)}+\|y_\Omega\|_{L^2(\Omega)}
\lesssim
\|f\|_{L^2(\Omega)}+\|y_\Omega\|_{L^2(\Omega)}.
\end{equation*}
These estimates, on the basis of \eqref{eq:step1_ii}, reveal that
\begin{equation}\label{eq:step1_iii}
\|\bar{u}-\tilde{u}\|_{L^2(\Omega)}
\lesssim
\|\nabla(\bar{p}_\T-\tilde{p})\|_{L^2(\Omega)}+\|\nabla(\bar{y}_\T-\tilde{y})\|_{L^2(\Omega)},
\end{equation}
with a hidden constant that is independent of continuous and discrete optimal variables but depends on the continuous problem data.

We now proceed to estimate $\|\nabla(\bar{y}_\T-\tilde{y})\|_{L^2(\Omega)}$ in \eqref{eq:step1_iii}. Invoke the auxiliary variable $\hat{y}\in H_0^1(\Omega)$, defined as the solution to \eqref{eq:hat_functions}, and the a posteriori error estimate \eqref{eq:estimate_state_hat_discrete_st}, to immediately arrive at
\begin{equation}\label{eq:step1_iv}
\|\nabla(\bar{y}_\T-\tilde{y})\|_{L^2(\Omega)} 
\lesssim
\|\nabla (\tilde{y}-\hat{y})\|_{L^2(\Omega)} + \mathcal{E}_{st,\T}.
\end{equation} 
To control $\|\nabla (\tilde{y}-\hat{y})\|_{L^2(\Omega)}$, we first observe that $\tilde{y}-\hat{y}\in H_0^1(\Omega)$ solves
\begin{equation*}
(\nabla (\tilde{y}-\hat{y}),\nabla v)_{L^2(\Omega)} + (\tilde{u}(\tilde{y}-\hat{y}),v)_{L^2(\Omega)}=(\hat{y}(\bar{u}_\T-\tilde{u}),v)_{L^2(\Omega)}\quad \forall v\in H_0^1(\Omega).
\end{equation*}
We can thus obtain, in view of a generalized H\"older's inequality, an stability bound for the problem that $\hat{y}$ solves, and definition \eqref{eq:error_estimator_control} the following estimate:
\[
\|\nabla(\tilde{y}-\hat{y})\|_{L^2(\Omega)} 
\lesssim
\|\nabla\hat{y}\|_{L^2(\Omega)}\|\bar{u}_\T-\tilde{u}\|_{L^2(\Omega)}
\lesssim 
\|f\|_{L^2(\Omega)}\|\bar{u}_\T-\tilde{u}\|_{L^2(\Omega)}
\lesssim
\mathcal{E}_{ct,\T}.
\]
Replacing this inequality into \eqref{eq:step1_iv} and the obtained one into \eqref{eq:step1_iii} yield
\begin{equation}\label{eq:step1_v}
\|\bar{u}-\tilde{u}\|_{L^2(\Omega)}
\lesssim
\|\nabla(\bar{p}_\T-\tilde{p})\|_{L^2(\Omega)} + \mathcal{E}_{ct,\T} + \mathcal{E}_{st,\T}.
\end{equation}

The rest of this step is dedicated to bound the term $\|\nabla(\bar{p}_\T-\tilde{p})\|_{L^2(\Omega)}$ in \eqref{eq:step1_v}. To accomplish this task, we first invoke the auxiliary variable $\hat{p}\in H_0^1(\Omega)$, defined as the solution to \eqref{eq:hat_functions_p} and the a posteriori error estimate \eqref{eq:estimate_state_hat_discrete_adj} to obtain
\begin{equation}\label{eq:step1_vi}
\|\nabla(\bar{p}_\T-\tilde{p})\|_{L^2(\Omega)} 
\lesssim
\mathcal{E}_{adj,\T} + \|\nabla(\tilde{p}-\hat{p})\|_{L^2(\Omega)}.
\end{equation}
To bound $ \|\nabla(\tilde{p}-\hat{p})\|_{L^2(\Omega)}$, we notice that
$\tilde{p}-\hat{p}\in H_0^1(\Omega)$ solves
\begin{multline*}
(\nabla w,\nabla (\tilde{p}-\hat{p}))_{L^2(\Omega)}+(\tilde{u}(\tilde{p}-\hat{p}),w)_{L^2(\Omega)}
\\
=
((\bar{u}_\T-\tilde{u})\hat{p},w)_{L^2(\Omega)} + (\tilde{y}-\bar{y}_\T,w)_{L^2(\Omega)}
\quad
\forall w\in H_0^1(\Omega).
\end{multline*}
This and the application of basic inequalities reveal the estimate
\begin{equation}\label{eq:step1_vi_half}
\|\nabla (\tilde{p}-\hat{p})\|_{L^2(\Omega)}
\lesssim
\|\nabla \hat{p}\|_{L^2(\Omega)}\|\bar{u}_\T-\tilde{u}\|_{L^2(\Omega)}+ \|\nabla (\tilde{y}-\bar{y}_\T)\|_{L^2(\Omega)}.
\end{equation}
Observe that 
\begin{equation}\label{eq:estimate_hat_p}
\|\nabla \hat{p}\|_{L^2(\Omega)}
\lesssim
\|\nabla \bar{y}_\T\|_{L^2(\Omega)} + \|y_\Omega\|_{L^2(\Omega)}
\lesssim 
\|f\|_{L^2(\Omega)} + \|y_\Omega\|_{L^2(\Omega)}.
\end{equation}
This estimate, combined with definition \eqref{eq:error_estimator_control}, and estimate $\|\nabla (\tilde{y}-\bar{y}_\T)\|_{L^2(\Omega)}\lesssim \mathcal{E}_{st,\T} + \mathcal{E}_{ct,\T}$ yield, on the basis of \eqref{eq:step1_vi_half}, the bound
\begin{equation}\label{eq:step1_vii}
\|\nabla (\tilde{p}-\hat{p})\|_{L^2(\Omega)}
\lesssim
\mathcal{E}_{ct,\T}+\mathcal{E}_{st,\T}.
\end{equation}
Replace estimate \eqref{eq:step1_vii} into \eqref{eq:step1_vi} and the obtained one into \eqref{eq:step1_v} to arrive at $\|\bar{u}-\tilde{u}\|_{L^2(\Omega)} \lesssim \mathcal{E}_{st,\T} + \mathcal{E}_{adj,\T} + \mathcal{E}_{ct,\T}$. Estimate \eqref{eq:step1_i} yields the desired bound
\begin{equation}\label{eq:step1_ix}
\|e_{\bar{u}}\|_{L^2(\Omega)}
\lesssim
\mathcal{E}_{st,\T} + \mathcal{E}_{adj,\T} + \mathcal{E}_{ct,\T}.
\end{equation}

%%%%%%%%%%%%%%%%%%%%%%%%%%%%%%%%%

\underline{Step 2.} The goal of this step is to control $\|\nabla e_{\bar{y}}\|_{L^2(\Omega)}$. Invoke the auxiliary varible $\hat{y}\in H_0^1(\Omega)$ and the a posteriori error estimate \eqref{eq:estimate_state_hat_discrete_st} to arrive at
\begin{equation}\label{eq:step2_i}
\|\nabla e_{\bar{y}}\|_{L^2(\Omega)} 
\lesssim
\|\nabla(\bar{y}-\hat{y})\|_{L^2(\Omega)}+\mathcal{E}_{st,\T}.
\end{equation}
To estimate $\|\nabla(\bar{y}-\hat{y})\|_{L^2(\Omega)}$, we first notice that $\bar{y}-\hat{y}\in H_0^1(\Omega)$ solves
\[
(\nabla(\bar{y}-\hat{y}),\nabla v)_{L^2(\Omega)}+(\bar{u}(\bar{y}-\hat{y}),v)_{L^2(\Omega)}=((\bar{u}_\T-\bar{u})\hat{y},v)_{L^2(\Omega)}\quad \forall v\in H_0^1(\Omega).
\]
Consequently, estimate \eqref{eq:step1_ix} yields
\begin{align*}
\|\nabla(\bar{y}-\hat{y})\|_{L^2(\Omega)} 
&\lesssim
\|\nabla \hat{y}\|_{L^2(\Omega)}\|\bar{u}_\T-\bar{u}\|_{L^2(\Omega)}\\
&\lesssim 
\|f\|_{L^2(\Omega)}\|e_{\bar{u}}\|_{L^2(\Omega)}
\lesssim \mathcal{E}_{st,\T} + \mathcal{E}_{adj,\T} + \mathcal{E}_{ct,\T}
.\nonumber
\end{align*}
We finally replace this bound into \eqref{eq:step2_i} to conclude that
\begin{equation}\label{eq:step2_ii}
\|\nabla e_{\bar{y}}\|_{L^2(\Omega)}
\lesssim
\mathcal{E}_{st,\T} + \mathcal{E}_{adj,\T} + \mathcal{E}_{ct,\T}.
\end{equation}

%%%%%%%%%%%%%%%%%%%%%%%%%%%%%%%%%%

\underline{Step 3.}  The objective now is to bound the term $\|\nabla e_{\bar{p}}\|_{L^2(\Omega)}$. To accomplish this task, we invoke $\hat{p}\in H_0^1(\Omega)$ and the bound \eqref{eq:estimate_state_hat_discrete_adj} to immediately arrive at
\begin{equation}\label{eq:step3_i}
\|\nabla e_{\bar{p}}\|_{L^2(\Omega)} 
\lesssim
\|\nabla(\bar{p}-\hat{p})\|_{L^2(\Omega)}+\mathcal{E}_{adj,\T}.
\end{equation}
To estimate $\|\nabla(\bar{p}-\hat{p})\|_{L^2(\Omega)}$, we observe that $\bar{p}-\hat{p}\in H_0^1(\Omega)$ solves
\[
\EO{(\nabla w,\nabla(\bar{p}-\hat{p}))_{L^2(\Omega)}}
+(\bar{u}(\bar{p}-\hat{p}),w)_{L^2(\Omega)}=((\bar{u}_\T-\bar{u})\hat{p},w)_{L^2(\Omega)} + (\bar{y}-\bar{y}_\T,w)_{L^2(\Omega)}
\]
for all $w\in H_0^1(\Omega)$. We thus invoke estimates \eqref{eq:estimate_hat_p}, \eqref{eq:step1_ix}, and \eqref{eq:step2_ii} to obtain
\begin{equation*}
\|\nabla(\bar{p}-\hat{p})\|_{L^2(\Omega)}
\lesssim
\|\nabla \hat{p}\|_{L^2(\Omega)}\|e_{\bar{u}}\|_{L^2(\Omega)}+\|e_{\bar{y}}\|_{L^2(\Omega)}\\
\lesssim \mathcal{E}_{st,\T} + \mathcal{E}_{adj,\T} + \mathcal{E}_{ct,\T}.
\end{equation*}
Replacing this bound into \eqref{eq:step3_i} we obtain 
\begin{equation}\label{eq:step3_ii}
\|\nabla e_{\bar{p}}\|_{L^2(\Omega)}\lesssim \mathcal{E}_{st,\T} + \mathcal{E}_{adj,\T} + \mathcal{E}_{ct,\T}.
\end{equation}

%%%%%%%%%%%%%%%%%%%%%%%%%%%%%%%%%

\underline{Step 4.} The desired estimate \eqref{eq:global_rel} follows from collecting the estimates \eqref{eq:step1_ix}, \eqref{eq:step2_ii}, and \eqref{eq:step3_ii}.
\qed

%%%%%%%%%%%%%%%%%%%%%%%%%%%%%%%%%%%%%%%%%%%%%%%%%%%%%%%%%%%%
%%%%%%%%%%%%%%%%%%%%%%%%%%%%%%%%%%%%%%%%%%%%%%%%%%%%%%%%%%%%
%%%%%%%%%%%%%%%%%%%%%%%%%%%%%%%%%%%%%%%%%%%%%%%%%%%%%%%%%%%%
%%%%%%%%%%%%%%%%%%%%%%%%%%%%%%%%%%%%%%%%%%%%%%%%%%%%%%%%%%%%

\subsection{Efficiency Analysis}\label{eq:sec:eff}

In this section, we derive local efficiency estimates for the local a posteriori error indicators $\mathcal{E}_{st,T}$ and $\mathcal{E}_{ad,T}$ and a global efficiency estimate for the a posteriori error estimator $\mathcal{E}_{ocp}$. To accomplish this task, we will proceed on the basis of standard residual estimation techniques \cite{MR3059294}. 

We begin our analysis by introducing the following notation: for an edge/face or triangle/tetrahedron $G$, we denote by $\mathcal{V}(G)$ the set of vertices of $G$. With this notation at hand, we introduce, for $T\in\mathscr{T}$ and $S\in\mathscr{S}$, the standard element and edge bubble functions \cite{MR3059294}, respectively, as 
\begin{equation}\label{eq:bubble_func}
\varphi^{}_{T}=
(d+1)^{(d+1)}\prod_{\textsc{v} \in \mathcal{V}(T)} \lambda^{}_{\textsc{v}},
\qquad
\varphi^{}_{S}=
d^{d} \prod_{\textsc{v} \in \mathcal{V}(S)}\lambda^{}_{\textsc{v}}|^{}_{T'}, 
\quad 
T' \subset \mathcal{N}_{S}.
\end{equation}
In these formulas, by $\lambda_\textsc{v}^{}$ we denote the barycentric coordinates of $T$. We recall that $\mathcal{N}_{S}$ denotes the patch composed of the two elements of $\mathscr{T}$ that share 
$S$.

The following identities are essential to perform an efficiency analysis. \FF{First,} since $\bar{y} \in H_0^1(\Omega)$ solves \eqref{eq:weak_state_eq} with $u=\bar{u}$, an elementwise integration by parts formula yields, for $v\in H_0^1(\Omega)$, the identity
\begin{multline}\label{eq:error_eq_state}
(\nabla e_{\bar{y}},\nabla v)_{L^2(\Omega)}+(e_{\bar{u}}\bar{y},v)_{L^2(\Omega)} = -(\bar{u}_\T e_{\bar{y}},v)_{L^2(\Omega)} +\sum_{S\in\Sides}(\llbracket\nabla \bar{y}_\T\cdot \boldsymbol{\nu}\rrbracket,v)_{L^2(S)}\\
+ \sum_{T\in\T}\left[(\mathscr{P}_{T}f-\bar{y}_\T\bar{u}_\T,v)_{L^2(T)} + (f-\mathscr{P}_{T}f,v)_{L^2(T)} \right].
\end{multline}
Second, since $\bar{p}$ solves \eqref{eq:adj_eq} with $y=\bar{y}$ and $u = \bar{u}$, similar arguments yield
\begin{multline}\label{eq:error_eq_adjoint}
(\nabla w,\nabla e_{\bar{p}})_{L^2(\Omega)}+(e_{\bar{u}}\bar{p},w)_{L^2(\Omega)} = -(\bar{u}_\T e_{\bar{p}},w)_{L^2(\Omega)}
+(e_{\bar{y}},w)_{L^2(\Omega)}
\\
+ \sum_{S\in\Sides}(\llbracket\nabla \bar{p}_\T\cdot \boldsymbol{\nu}\rrbracket,w)_{L^2(S)}
+ \sum_{T\in\T}  (\bar{y}_\T-\mathscr{P}_{T}y_\Omega-\bar{u}_\T\bar{p}_\T,w)_{L^2(T)} 
\\
+\sum_{T\in\T}  (\mathscr{P}_{T}y_\Omega-y_\Omega,w)_{L^2(\FF{T})}
\quad
\forall w\in H_0^1(\Omega).
\end{multline}
In \eqref{eq:error_eq_state} and \eqref{eq:error_eq_adjoint}, $\mathscr{P}_{T}$ denotes the $L^2$-projection operator onto piecewise constant functions over $T\in\T$.

As a final ingredient, we introduce, for $v \in L^2(\Omega)$ and $\mathcal{M}\subset \T$,
\begin{equation*}
\mathrm{osc}_{\T}(v;\mathcal{M}):= \left(\sum_{T\in\mathcal{M}}h_T^2\|v-\mathscr{P}_{T}v\|_{L^2(T)}^2\right)^\frac{1}{2}.
\end{equation*}

We are now ready to prove the local efficiency of $\mathcal{E}_{st,T}$, defined in \eqref{eq:indicators_state_eq}.

\begin{theorem}[local efficiency of $\mathcal{E}_{st\FF{,T}}$]\label{thm:local_eff_st} 
Let $\bar{u} \in \mathbb{U}_{ad}$ be a local solution to \eqref{eq:minimize_weak}--\eqref{eq:weak_state_eq}.
% satisfying the \FF{first order necessary condition \eqref{eq:var_ineq}}. 
Let $\bar u_{\T}$ be a local minimum of the fully discrete optimal control problem with $\bar{y}_{\T}$ and $\bar{p}_\T$ being the corresponding state and adjoint state, respectively. Then, for $T\in\T$, the local error indicator $\mathcal{E}_{st,T}$ satisfies the  bound
\begin{equation}\label{eq:local_eff_st}
\mathcal{E}_{st,T}
\lesssim
\|\nabla e_{\bar{y}}\|_{L^2(\mathcal{N}_{T})} + \|e_{\bar{u}}\|_{L^2(\mathcal{N}_{T})} + h_{T}\|e_{\bar{y}}\|_{L^2(\mathcal{N}_T)} + \mathrm{osc}_{\T}(f;\mathcal{N}_T),
\end{equation}
where $\mathcal{N}_{T}$ is defined in \eqref{def:patch}. The hidden constant is independent of  continuous and discrete optimal variables, the size of the elements in $\T$, and $\#\T$.
\end{theorem}
\proof
We proceed in two steps and estimate each term in the definition of the local error indicator $\mathcal{E}_{st,T}$, given in \eqref{eq:indicators_state_eq}, separately.

\underline{Step 1.} Let $T\in\T$. We bound $h_{T}^2\|f-\bar{y}_\T\bar{u}_\T\|_{L^2(T)}^2$ in \eqref{eq:indicators_state_eq}. To accomplish this task, we begin with a simple application of a triangle inequality to write
\begin{equation}\label{eq:triangle_st_eff}
h_{T}^2\|f-\bar{y}_\T\bar{u}_\T\|_{L^2(T)}^2
\leq 
2h_{T}^2\|\mathscr{P}_{T}f-\bar{y}_\T\bar{u}_\T\|_{L^2(T)}^2 + 2\mathrm{osc}_{\T}(f;T)^2.
\end{equation}
It thus suffices to bound the term $h_T\|\mathscr{P}_{T}f-\bar{y}_\T\bar{u}_\T\|_{L^2(T)}$. To do this, we set $v=\varphi_{T}(\mathscr{P}_{T}f-\bar{y}_\T\bar{u}_\T)$ in \eqref{eq:error_eq_state}, where $\varphi_T$ denotes the element bubble function introduced in \eqref{eq:bubble_func}. Standard properties of $\varphi_T$ and inverse inequalities yield
\begin{multline*}
\|\mathscr{P}_{T}f-\bar{y}_\T\bar{u}_\T\|_{L^2(T)}^2
\lesssim
\left( h_{T}^{-1}\|\nabla e_{\bar{y}}\|_{L^2(T)}+h_{T}^{-1}\|\bar{y}\|_{L^{d}(T)}\|e_{\bar{u}}\|_{L^2(T)}\right.\\
\left.+\|\bar{u}_\T\|_{L^\infty(T)}\|e_{\bar{y}}\|_{L^2(T)}+\|f-\mathscr{P}_{T}f\|_{L^2(T)}\right)\|\mathscr{P}_{T}f-\bar{y}_\T\bar{u}_\T\|_{L^2(T)}.
\end{multline*}
We notice that, since $d \in \{2,3\}$, $H_0^1(\Omega)\hookrightarrow L^d(\Omega)$. Consequently,
\begin{equation}\label{eq:embedding_bary}
\|\bar{y}\|_{L^{d}(T)} \leq \|\bar{y}\|_{L^{d}(\Omega)} \lesssim \|\nabla \bar{y}\|_{L^2(\Omega)} \lesssim \|f\|_{L^2(\Omega)}.
\end{equation}
This bound combined with the fact that $\bar{u}_\T\in \mathbb{U}_{ad}$ yield
\begin{equation*}
h_{T}\|\mathscr{P}_{T}f-\bar{y}_\T\bar{u}_\T\|_{L^2(T)}
\lesssim
\|\nabla e_{\bar{y}}\|_{L^2(T)} + \|e_{\bar{u}}\|_{L^2(T)} +h_{T}\|e_{\bar{y}}\|_{L^2(T)}+\mathrm{osc}_{\T}(f;T).
\end{equation*}
Replace this bound into \eqref{eq:triangle_st_eff} to obtain the desired one for $h_{T}^2\|f-\bar{y}_\T\bar{u}_\T\|_{L^2(T)}^2$.

\underline{Step 2.} Let $T\in\T$ and $S\in\mathscr{S}_{T}$. We now bound $h_T^{\frac{1}{2}}\|\llbracket\nabla \bar{y}_\T\cdot \mathbf{\boldsymbol{\nu}} \rrbracket\|_{L^2(S)}$ in \eqref{eq:indicators_state_eq}. As a first step, we set $v=\varphi_{S}\llbracket\nabla \bar{y}_\T\cdot \boldsymbol{\nu}\rrbracket$ in identity \eqref{eq:error_eq_state}. Here, $\varphi_{S}$ denotes the edge bubble function introduced in  \eqref{eq:bubble_func}. We thus invoke standard bubble functions arguments and inverse inequalities to arrive at
\begin{multline*}
\|\llbracket\nabla \bar{y}_\T\cdot \boldsymbol{\nu}\rrbracket\|_{L^2(S)}^2
\lesssim
\sum_{T'\in\mathcal{N}_{T}}\big( h_{T'}^{-1}\|\nabla e_{\bar{y}}\|_{L^2(T')}
+
h_{T'}^{-1} \|\bar{y}\|_{L^{d}(T')}\|e_{\bar{u}}\|_{L^2(T')}
\\
+
\|\bar{u}_\T\|_{L^\infty(T')}\|e_{\bar{y}}\|_{L^2(T')} 
+ \|f-\mathscr{P}_{T'}f\|_{L^2(T')} 
\\
+ \|\mathscr{P}_{T}f-\bar{y}_\T\bar{u}_\T\|_{L^2(T')}\big)h_{T}^{\frac{1}{2}}\|\llbracket\nabla \bar{y}_\T\cdot \boldsymbol{\nu}\rrbracket\|_{L^2(S)}.
\end{multline*}
In view of the derived estimate for $h_T\|\mathscr{P}_{T}f-\bar{y}_\T\bar{u}_\T\|_{L^2(T)}$, \eqref{eq:embedding_bary}, and the fact that $\bar{u}_\T\in \mathbb{U}_{ad}$, we immediately conclude that
\begin{multline*}
h_{T}^{\frac{1}{2}}\|\llbracket\nabla \bar{y}_\T\cdot \boldsymbol{\nu}\rrbracket\|_{L^2(S)}
\lesssim
\sum_{T'\in\mathcal{N}_{T}}\big( \|\nabla e_{\bar{y}}\|_{L^2(T')}+\|e_{\bar{u}}\|_{L^2(T')}\\
+h_{T}\|e_{\bar{y}}\|_{L^2(T')}+\mathrm{osc}_{\T}(f;T')\big).
\end{multline*}

A collection of the bounds derived in Steps 1 and 2 yield \eqref{eq:local_eff_st}. This concludes the proof.
\qed

We now continue with the study of local efficiency properties for the indicator $\mathcal{E}_{adj,T}$, which is defined in \eqref{eq:error_indicator_adjoint_eq}.

\begin{theorem}[local efficiency of $\mathcal{E}_{adj\FF{,T}}$]\label{thm:local_eff_adj}
In the framework of Theorem \ref{thm:local_eff_st}, we have, for $T\in\T$, the local estimate
\begin{multline}\label{eq:local_eff_adj}
\mathcal{E}_{adj,T}
\lesssim
\|\nabla e_{\bar{p}}\|_{L^2(\mathcal{N}_{T})} +  \|e_{\bar{u}}\|_{L^2(\mathcal{N}_{T})}  \\
+ h_{T}\|e_{\bar{p}}\|_{L^2(\mathcal{N}_T)}
+ h_{T}\|e_{\bar{y}}\|_{L^2(\mathcal{N}_{T})} + \mathrm{osc}_{\T}(y_\Omega;\mathcal{N}_T),
\end{multline}
where $\mathcal{N}_{T}$ is defined in \eqref{def:patch}. The hidden constant is independent of continuous and discrete optimal variables, the size of the elements in $\T$, and $\#\T$.
\end{theorem}
\proof
The proof relies on utilizing identity \eqref{eq:error_eq_adjoint} and similar arguments to the ones elaborated within the proof of Theorem \ref{thm:local_eff_st}. For brevity, we skip details.
\qed

The results of Theorems \ref{thm:local_eff_st} and \ref{thm:local_eff_adj} yield the following global efficiency estimate for $\mathcal{E}_{ocp,\T}$. 

\begin{theorem}[global efficiency of $\mathcal{E}_{ocp\FF{,\T}}$]\label{thm:global_eff}
In the framework of Theorem \ref{thm:local_eff_st}, we have the global estimate
\begin{equation*}
\mathcal{E}_{ocp,\T}
\lesssim 
\|\nabla e_{\bar{p}}\|_{L^2(\Omega)}+\|\nabla e_{\bar{y}}\|_{L^2(\Omega)}+\|e_{\bar{u}}\|_{L^2(\Omega)} + \mathrm{osc}_{\T}(f;\T)+ \mathrm{osc}_{\T}(y_{\Omega};\T).
\end{equation*}
The hidden constant is independent of continuous and discrete optimal variables, the size of the elements in $\T$, and $\#\T$.
\end{theorem}
\proof
We begin by invoking the definition of the error estimator $\mathcal{E}_{st,\T}$, given in \eqref{eq:error_estimator_state_eq}, and the local efficiency estimate \eqref{eq:local_eff_st} to arrive at
\begin{equation}\label{eq:global_estimate_st}
\mathcal{E}_{st,\T}
\lesssim 
\|\nabla e_{\bar{y}}\|_{L^2(\Omega)} + \|e_{\bar{u}}\|_{L^2(\Omega)} + \text{diam}(\Omega)\|e_{\bar{y}}\|_{L^2(\Omega)} + \mathrm{osc}_{\T}(f,\T).
\end{equation}
On the other hand, the definition of the error estimator $\mathcal{E}_{adj,\T}$, given in \eqref{eq:error_estimator_adjoint_eq}, and the efficiency estimate \eqref{eq:local_eff_adj} yield the bound
\begin{multline}\label{eq:global_estimate_adj}
\mathcal{E}_{adj,\T}
\lesssim 
\|\nabla e_{\bar{p}}\|_{L^2(\Omega)} + \|e_{\bar{u}}\|_{L^2(\Omega)} +  \text{diam}(\Omega)\|e_{\bar{p}}\|_{L^2(\Omega)} \\
+ \text{diam}(\Omega)\|e_{\bar{y}}\|_{L^2(\Omega)} + \mathrm{osc}_{\T}(y_\Omega,\T).
\end{multline}
It thus suffices to bound the estimator $\mathcal{E}_{ct,\T}$. In view of \eqref{eq:error_estimator_control}, a trivial application of a triangle inequality yields
\begin{align*}
\mathcal{E}_{ct,\T} 
&\leq \|\tilde{u}-\bar{u}\|_{L^2(\Omega)}+\|e_{\bar{u}}\|_{L^2(\Omega)}\\
&= \|\Pi_{[\texttt{a},\texttt{b}]}(\alpha^{-1}\bar{p}_\T\bar{y}_\T)-\Pi_{[\texttt{a},\texttt{b}]}(\alpha^{-1}\bar{p}\bar{y})\|_{L^2(\Omega)}+\|e_{\bar{u}}\|_{L^2(\Omega)},
\end{align*}
where $\Pi_{[\texttt{a},\texttt{b}]}$ is defined in \eqref{def:projector_pi}. This bound, the Lipschitz property of $\Pi_{[\texttt{a},\texttt{b}]}$, the Cauchy--Schwarz inequality, and the embedding $H_0^1(\Omega)\hookrightarrow L^4(\Omega)$ yield
\begin{equation*}
\mathcal{E}_{ct,\T}
\leq
\alpha^{-1}\|\nabla e_{\bar{p}}\|_{L^2(\Omega)}\|\nabla \bar{y}_\T\|_{L^2(\Omega)}  +\alpha^{-1}\|\nabla\bar{p}\|_{L^2(\Omega)}\|\nabla e_{\bar{y}}\|_{L^2(\Omega)}+\|e_{\bar{u}}\|_{L^2(\Omega)}.
\end{equation*}
Observe that $\|\nabla \bar{y}_\T\|_{L^2(\Omega)}\lesssim \|f\|_{L^2(\Omega)}$ and $\|\nabla\bar{p}\|_{L^2(\Omega)}\lesssim \|f\|_{L^2(\Omega)} +\|y_\Omega\|_{L^2(\Omega)}$. Consequently,
\begin{equation}\label{eq:global_estimate_ct}
\mathcal{E}_{ct,\T} 
\lesssim
\|\nabla e_{\bar{p}}\|_{L^2(\Omega)}+\|\nabla e_{\bar{y}}\|_{L^2(\Omega)}+\|e_{\bar{u}}\|_{L^2(\Omega)},
\end{equation}
with a hidden constant that is independent of continuous and discrete optimal variables but depends on the continuous problem data.

The proof concludes by gathering estimates \eqref{eq:global_estimate_st}, \eqref{eq:global_estimate_adj}, and \eqref{eq:global_estimate_ct}, upon utilizing a Poincar\'e inequality.
\qed

%%%%%%%%%%%%%%%%%%%%%%%%%%%%%%%%%%%%%%%%%%%%%%%%%%%%%%%%%%%%
%%%%%%%%%%%%%%%%%%%%%%%%%%%%%%%%%%%%%%%%%%%%%%%%%%%%%%%%%%%%
%%%%%%%%%%%%%%%%%%%%%%%%%%%%%%%%%%%%%%%%%%%%%%%%%%%%%%%%%%%%
%%%%%%%%%%%%%%%%%%%%%%%%%%%%%%%%%%%%%%%%%%%%%%%%%%%%%%%%%%%%
%%%%%%%%%%%%%%%%%%%%%%%%%%%%%%%%%%%%%%%%%%%%%%%%%%%%%%%%%%%%
%%%%%%%%%%%%%%%%%%%%%%%%%%%%%%%%%%%%%%%%%%%%%%%%%%%%%%%%%%%%
%%%%%%%%%%%%%%%%%%%%%%%%%%%%%%%%%%%%%%%%%%%%%%%%%%%%%%%%%%%%
%%%%%%%%%%%%%%%%%%%%%%%%%%%%%%%%%%%%%%%%%%%%%%%%%%%%%%%%%%%%

\section{A Posteriori Error Analysis: the Semi-discrete Scheme}\label{sec:a_posteriori_semi}

In this section, we design and analyze an a posteriori error estimator for the semi-discrete scheme of section \ref{sec:semi_discrete_scheme}. In contrast to the estimator devised in section \ref{sec:a_posteriori_fully}, the estimator is now \FF{formed by} only two contributions: one related to the discretization of the state equation and another one associated to the discretization of the adjoint equation.

\subsection{Global Reliability Analysis}\label{sec:global_var}

The goal of this section is to design an a posteriori error estimator and derive an upper bound for the corresponding total error in terms of the devised error estimator. As in section \ref{sec:global_fully}, the aforementioned upper bound will be obtained on the basis of estimates on the error between solutions to the semi-discrete optimal control problem \eqref{eq:semi_discrete_state_equation}--\eqref{eq:semi_discrete_adjoint_equation} and suitable auxiliary variables.

The first auxiliary variable is $\hat y$ and is defined as follows:
\begin{equation}\label{eq:hat_functions_var}
\hat y \in H_{0}^{1}(\Omega):
\quad
(\nabla\hat{y},\nabla v)_{L^2(\Omega)}+(\bar{\mathsf{u}}\hat{y},v)_{L^2(\Omega)} =  (f,v)^{}_{L^2(\Omega)} \quad \forall v \in H_0^1(\Omega).
\end{equation}
With this variable at hand, we define, for $T\in\T$, the local error indicators and the corresponding a posteriori error estimator, respectively, by
\begin{equation*}
\mathsf{E}_{st,T}^2:=h_T^2\|f-\bar{\mathsf{u}}\bar{y}_\T\|_{L^2(T)}^2+h_T\|\llbracket\nabla \bar{y}_\T\cdot\boldsymbol{\nu}\rrbracket\|_{L^2(\partial T\setminus\partial\Omega)}^2, \quad \mathsf{E}_{st,\T}^2:=\sum_{T\in\T}\mathsf{E}_{st,T}^2.
\end{equation*}
We notice that, since $\bar{y}_{\T}\in\mathbb{V}(\T)$, solution to \eqref{eq:semi_discrete_state_equation} with $\mathsf{u}=\bar{\mathsf{u}}$, can be seen as a finite element approximation of $\hat{y}$, an application of Theorem \ref{thm:global_reli_weak} yields
\begin{equation}\label{eq:estimate_state_hat_discrete_st_var}
\|\nabla (\hat{y}-\bar{y}_\T)\|_{L^2(\Omega)}
\lesssim \mathsf{E}_{st,\T}.
\end{equation}

The second variable is $\hat p \in H_{0}^{1}(\Omega)$ and is defined as the solution to
\begin{equation}\label{eq:hat_functions_p_var}
(\nabla w,\nabla \hat{p})_{L^2(\Omega)}+\left(\bar{\mathsf{u}}\hat{p},w\right)_{L^2(\Omega)} =  (\bar{y}_\T-y^{}_{\Omega},w)^{}_{L^2(\Omega)} \quad \forall w\in H_0^1(\Omega).
\end{equation}
Define, for $T\in\T$, the local error indicators
\begin{equation*}
\mathsf{E}_{adj,T}^2:=h_T^2\|\bar{y}_\T-y_\Omega-\bar{\mathsf{u}}\bar{p}_{\T}\|_{L^2(T)}^2+h_T\|\llbracket\nabla \bar{p}_\T\cdot\boldsymbol{\nu}\rrbracket\|_{L^2(\partial T\setminus\partial\Omega)}^2,
\end{equation*}
and the a posteriori error estimator
\begin{equation*}
\mathsf{E}_{adj,\T}^2:=\sum_{T\in\T}\mathsf{E}_{adj,T}^2.
\end{equation*}
\EO{Since $\bar{p}_{\T}\in\mathbb{V}(\T)$
% , solution to \eqref{eq:semi_discrete_adjoint_equation} with $\mathsf{u}=\bar{\mathsf{u}}$, 
can} be seen as the finite element approximation of $\hat{p}$ within $\mathbb{V}(\T)$, Theorem \ref{thm:global_reli_weak} yields
\begin{equation}\label{eq:estimate_state_hat_discrete_adj_var}
\|\nabla (\hat{p}-\bar{p}_\T)\|_{L^2(\Omega)}
\lesssim \mathsf{E}_{adj,\T}.
\end{equation}

In order to present the following reliability result, we introduce the error $\mathsf{e}_{\bar{u}}:= \bar{u} - \bar{\mathsf{u}}$ and the a posteriori error estimator $\mathsf{E}_{ocp,\T}^2:=\mathsf{E}_{st,\T}^2 + \mathsf{E}_{adj,\T}^2$.

\begin{theorem}[global reliability]\label{thm:global_rel_var}
Let $\bar u \in \mathbb{U}_{ad}$ be a local solution to \eqref{eq:minimize_weak}--\eqref{eq:weak_state_eq} satisfying the sufficient second order condition \eqref{eq:second_order_cond}. Let $\bar{\mathsf{u}}$ be a local minimum of the semi-discrete optimal control problem with $\bar{y}_{\T}$ and $\bar{p}_\T$ being the corresponding state and adjoint state, respectively. If \FF{$\bar{y}_{\T}$ and $\bar{p}_\T$ satisfy, on the mesh $\T$, the bound \eqref{eq:assumption_mesh}}, then
\begin{equation}\label{eq:global_rel_var}
\|\nabla{e}_{\bar{y}}\|^2_{L^2(\Omega)} + \|\nabla{e}_{\bar{p}}\|^2_{L^2(\Omega)} + \| \mathsf{e}_{\bar{u}} \|^2_{L^2(\Omega)}
\lesssim 
\mathsf{E}_{ocp,\T}^2, 
\end{equation}
with a hidden constant that is independent of continuous and discrete optimal variables, the size of the elements in $\T$, and $\#\T$.
\end{theorem}
\proof
We immediately notice that, under the particular setting inherited by the semi-discrete scheme, the auxiliary variables that we have devised to perform our analysis satisfy that $\bar{\mathsf{u}}=\tilde{u}$, where $\tilde{u}$ is defined in \eqref{def:control_tilde}. We can thus immediately conclude that $\hat{y}=\tilde{y}$, with $\hat{y}$ and $\tilde{y}$ being defined as the unique solutions to \eqref{eq:hat_functions_var} and \eqref{eq:tilde_y}, respectively. Therefore, invoking \eqref{eq:step1_iii}, \eqref{eq:step1_iv}, \eqref{eq:step1_vi}, and \eqref{eq:step1_vi_half} we conclude that 
\begin{equation*}
\|\mathsf{e}_{\bar{u}}\|_{L^2(\Omega)} = \|\bar{u}-\tilde{u}\|_{L^2(\Omega)} 
\lesssim 
\|\nabla(\bar{p}_\T - \tilde{p})\|_{L^2(\Omega)} + \mathsf{E}_{st,\T}
\lesssim
\mathsf{E}_{st,\T}
+ \mathsf{E}_{adj,\T}.
\end{equation*}
The \FF{estimation} of the terms $\|\nabla{e}_{\bar{y}}\|_{L^2(\Omega)}$ and $\|\nabla{e}_{\bar{p}}\|_{L^2(\Omega)}$ follow by utilizing the bound $\|\mathsf{e}_{\bar{u}}\|_{L^2(\Omega)}\lesssim \mathsf{E}_{ocp,\T}$ and similar arguments to the ones developed in the proof of Theorem \ref{thm:global_rel}. For brevity, we skip the details. 
\qed

%%%%%%%%%%%%%%%%%%%%%%%%%%%%%%%%%%%%%%%%%%%%%%%%%%%%%%%%%%%%
%%%%%%%%%%%%%%%%%%%%%%%%%%%%%%%%%%%%%%%%%%%%%%%%%%%%%%%%%%%%
%%%%%%%%%%%%%%%%%%%%%%%%%%%%%%%%%%%%%%%%%%%%%%%%%%%%%%%%%%%%
%%%%%%%%%%%%%%%%%%%%%%%%%%%%%%%%%%%%%%%%%%%%%%%%%%%%%%%%%%%%

\subsection{Efficiency Analysis}

We begin the section by defining, for $T\in\T$, the local indicator
\begin{equation}\label{def:indicator_ocp_variational}
\mathsf{E}_{ocp,T}^2:= \mathsf{E}_{st,T}^{2} + \mathsf{E}_{adj,T}^{2}.
\end{equation}
The estimates obtained in Theorems \ref{thm:local_eff_st} and \ref{thm:local_eff_adj} can also be obtained within the setting of the variational discretization approach.

\begin{theorem}[local estimates for $\mathsf{E}_{ocp,T}$]
\label{thm:global_eff_var}  
Let $\bar u \in \mathbb{U}_{ad}$ be a local solution to \eqref{eq:minimize_weak}--\eqref{eq:weak_state_eq}.
% satisfying the \FF{first order necessary condition \eqref{eq:var_ineq}}. 
Let $\bar{\mathsf{u}}$ be a local minimum of the semi-discrete optimal control problem with $\bar{y}_{\T}$ and $\bar{p}_\T$ being the corresponding state and adjoint state, respectively. Then, for $T\in\T$, the local error indicator $\mathsf{E}_{ocp,T}$ satisfies
\begin{multline*}
\mathsf{E}_{ocp,T} 
\lesssim
\|\nabla e_{\bar{y}}\|_{L^2(\mathcal{N}_{T})} + \|\nabla e_{\bar{p}}\|_{L^2(\mathcal{N}_{T})} +  \|\mathsf{e}_{\bar{u}}\|_{L^2(\mathcal{N}_{T})} + h_{T}\|e_{\bar{p}}\|_{L^2(\mathcal{N}_T)} \\+ h_{T}\|e_{\bar{y}}\|_{L^2(\mathcal{N}_{T})} + \mathrm{osc}_{\T}(f;\mathcal{N}_T) + \mathrm{osc}_{\T}(y_\Omega;\mathcal{N}_T),
\end{multline*}
where $\mathcal{N}_T$ is defined as in \eqref{def:patch}. The hidden constant is independent of continuous and discrete optimal variables, the size of the elements in $\T$, and $\#\T$.
\end{theorem}
\proof
The desired result follows by utilizing similar arguments to the ones that yield estimates \eqref{eq:local_eff_st} and \eqref{eq:local_eff_adj}. For brevity, we skip details.
\qed

%%%%%%%%%%%%%%%%%%%%%%%%%%%%%%%%%%%%%%%%%%%%%%%%%%%%%%%%%%%%%%%%%%%
%%%%%%%%%%%%%%%%%%%%%%%%%%%%%%%%%%%%%%%%%%%%%%%%%%%%%%%%%%%%%%%%%%%
%%%%%%%%%%%%%%%%%%%%%%%%%%%%%%%%%%%%%%%%%%%%%%%%%%%%%%%%%%%%%%%%%%%
%%%%%%%%%%%%%%%%%%%%%%%%%%%%%%%%%%%%%%%%%%%%%%%%%%%%%%%%%%%%%%%%%%%
%%%%%%%%%%%%%%%%%%%%%%%%%%%%%%%%%%%%%%%%%%%%%%%%%%%%%%%%%%%%%%%%%%%
%%%%%%%%%%%%%%%%%%%%%%%%%%%%%%%%%%%%%%%%%%%%%%%%%%%%%%%%%%%%%%%%%%%
%%%%%%%%%%%%%%%%%%%%%%%%%%%%%%%%%%%%%%%%%%%%%%%%%%%%%%%%%%%%%%%%%%%
%%%%%%%%%%%%%%%%%%%%%%%%%%%%%%%%%%%%%%%%%%%%%%%%%%%%%%%%%%%%%%%%%%%

\section{Numerical Examples}\label{sec:numerical_ex}

In this section, we conduct a series of numerical experiments that illustrate the performance of the devised a posteriori error estimators $\mathcal{E}_{ocp,\T}$ and $\mathsf{E}_{ocp,\T}$ when used to drive suitable AFEMs schemes based on the fully and semi-discrete schemes proposed in sections \ref{sec:fully_discrete_scheme} and \ref{sec:semi_discrete_scheme}, respectively.

\subsection{Implementation Details}

The numerical examples that we \FF{shall present in what follows} have been carried out with the help of a code that we implemented using \texttt{C++}. Global linear systems were solved using the multifrontal massively parallel sparse direct solver (MUMPS) \cite{MUMPS1,MUMPS2}. \FF{We have used a quadrature formula to compute} the right-hand sides, the approximation errors, and the error indicators\FF{; the quadrature formula being} exact for polynomials of degree nineteen $(19)$ for two dimensional domains and degree fourteen $(14)$ for three dimensional domains. 

In what follows, we discuss some pertinent implementation details that are particular for each discretization technique.

\textbf{The fully discrete scheme}: For a given partition $\mathscr{T}$, we seek a discrete solution $(\bar{y}_\mathscr{T},\bar{p}^{}_\mathscr{T},\bar{u}_\mathscr{T})\in \mathbb{V}(\mathscr{T}) \times \mathbb{V}(\mathscr{T}) \times \mathbb{U}_{ad}(\T)$ that solves the corresponding optimality system. This nonlinear system is solved on the basis of an adaptation of the semi--smooth Newton method described in \cite[Appendix A.1]{MR2971171}. \FF{We have exactly assembled all matrices involved in the left-hand side of the resulting linear system}. The total number of degrees of freedom (DOFs) is $\mathsf{Ndof}=2\dim(\mathbb{V}(\T)) + \dim(\mathbb{U}(\T))$. \EO{We measure the error within the norm} $\|e\|_{\Omega}:=[\|\nabla{e}_{\bar{y}}\|^2_{L^2(\Omega)} + \|\nabla{e}_{\bar{p}}\|^2_{L^2(\Omega)} + \| e_{\bar{u}} \|^2_{L^2(\Omega)}]^{\frac{1}{2}}$. Finally, we introduce the effectivity index $\mathcal{I}_{eff}:= \mathcal{E}_{ocp,\T}/\|e\|_{\Omega}$.

\textbf{The semi-discrete scheme}: For a given partition $\mathscr{T}$, we seek a solution $(\bar{y}_\mathscr{T},\bar{p}^{}_\mathscr{T})\in \mathbb{V}(\mathscr{T}) \times \mathbb{V}(\mathscr{T}) $ that solves the corresponding optimality system. This system is also solved by using an adaptation of the semi--smooth Newton method described in \cite[Appendix A.1]{MR2971171}.  

\EO{The following comments regarding the implementation of the variational discretization approach are of importance. In order to properly implement such a scheme, the assembling and exact computation of $(\bar{\mathsf{u}}\bar{y}_{\T},v_{\mathscr{T}})_{L^{2}(\Omega)}$ and $(\bar{\mathsf{u}}\bar{p}_{\T},w_{\mathscr{T}})_{L^{2}(\Omega)}$ are required. In particular, the exact integration of such terms on the simplices $T\in\T$ where the control $\bar{\mathsf{u}}$ exhibits kinks is necessary. Let us now describe an alternative to perform such a computation:
% \FF{An} alternative to accomplish this task is as follows: 
First, recognize the simplices $T \in \T$ which are such that the control $\bar{\mathsf{u}}$ have kinks. Second, recognize the regions of such simplices where the control variable is inactive/active. The following difficulty thus appears: 
% We immediately} notice that, 
since we are using continuous piecewise polynomials of degree one to approximate $\bar{y}_{\T}$ and $\bar{p}_{\T}$, these regions have, in general, curved boundaries; see \cite[Remark 5.19]{MR2536007}. The third step is the computation of $(\bar{\mathsf{u}}\bar{y}_{\T},v_{\mathscr{T}})_{L^{2}(\Omega)}$ and $(\bar{\mathsf{u}}\bar{p}_{\T},w_{\mathscr{T}})_{L^{2}(\Omega)}$ by partitioning the integrals on the regions where the control is inactive/active. Since the computational implementation of the second and third steps is far from being simple, our implementation relies on computing the terms $(\bar{\mathsf{u}}\bar{y}_{\T},v_{\mathscr{T}})_{L^{2}(\Omega)}$ and $(\bar{\mathsf{u}}\bar{p}_{\T},w_{\mathscr{T}})_{L^{2}(\Omega)}$ with the help of a quadrature formula. We emphasize that this numerical implementation leads to an \emph{approximated version} of the variational discretization approach which is capable of delivering, within an adaptive loop, optimal experimental rates of convergence for all the involved variables.}
 
The total number of DOFs for the semi-discrete scheme corresponds to $\mathsf{Ndof}=2\dim(\mathbb{V}(\T))$. To measure the corresponding approximation error, we use $\|\mathsf{e}\|_{\Omega}=[\|\nabla{e}_{\bar{y}}\|^2_{L^2(\Omega)} + \|\nabla{e}_{\bar{p}}\|^2_{L^2(\Omega)} + \| \mathsf{e}_{\bar{u}} \|^2_{L^2(\Omega)}]^{\frac{1}{2}}$. As a final ingredient, we introduce the effectivity index $\mathsf{I}_{eff}:= \mathsf{E}_{ocp,\T}/\|\mathsf{e}\|_{\Omega}$. 

Once the discrete solution is obtained, we compute, for $T\in \T$, the error indicator $\mathcal{E}_{ocp,T}$, defined by
\begin{equation}\label{def:indicator_ocp}
\mathcal{E}_{ocp,T}^2:= \mathcal{E}_{st,T}^2+\mathcal{E}_{adj,T}^2+\mathcal{E}_{ct,T}^2,
\end{equation}
\EO{or} the indicator $\mathsf{E}_{ocp,T}$, defined in \eqref{def:indicator_ocp_variational}, to drive the adaptive mesh refinement procedure described in Algorithm \ref{Algorithm}. 

\begin{algorithm}[ht]
\caption{\textbf{Adaptive algorithm}}
\label{Algorithm}
\textbf{Input:} Initial mesh $\mathscr{T}_{0}$, desired state $y_\Omega$, external source $f$, constraints $\texttt{a}$ and $\texttt{b}$, and regularization parameter $\alpha$;
\\
\textbf{Set:} $i=0$.
\\
\textbf{Newton strategy:}
\\
$\boldsymbol{1}$: Choose initial guesses $y_{\T_i}^0,p_{\T_i}^0 \in \mathbb{V}(\T_i)$ (and $u_{\T_i}^0\in \mathbb{U}(\T_i)$ when the fully discrete scheme is considered); 
\\
$\boldsymbol{2}$ (Fully discrete solution technique): Compute $[\bar{y}_{\T_i},\bar{p}_{\T_i},\bar{u}_{\T_i}] = \text{Semi-Smooth}[\T_{i},y_{\T_i}^0,p_{\T_i}^0,u_{\T_i}^0,y_\Omega,f,\texttt{a},\texttt{b},\alpha]$, which implements an adaptation of the semi--smooth Newton method described in \cite[Appendix A.1]{MR2971171};
\\
$\boldsymbol{2}$ (Semi-discrete solution technique): Compute $[\bar{y}_{\T_i},\bar{p}_{\T_i}] = \text{Semi-Smooth}[\T_{i},y_{\T_i}^0,p_{\T_i}^0,y_\Omega,f,\texttt{a},\texttt{b},\alpha]$, which implements an adaptation of the semi--smooth Newton method described in \cite[Appendix A.1]{MR2971171};
\\
\textbf{Adaptive loop:}
\\
$\boldsymbol{3}$: For each $T \in \mathscr{T}_i$, compute the local error indicator $\mathcal{E}_{ocp,T}$ ($\mathsf{E}_{ocp,T}$) defined in \eqref{def:indicator_ocp} (\eqref{def:indicator_ocp_variational});
\\
$\boldsymbol{4}$: Mark an element $T \in \T_i$ for refinement if $\mathcal{E}_{ocp,T}^{2}> \frac{1}{2}\max_{T'\in \mathscr{T}_i}\mathcal{E}_{ocp,T'}^{2}$ ($\mathsf{E}_{ocp,T}^{2}> \frac{1}{2}\max_{T'\in \mathscr{T}_i}\mathsf{E}_{ocp,T'}^{2}$);
\\
$\boldsymbol{5}$: From step $\boldsymbol{4}$, construct a new mesh $\mathscr{T}_{i+1}$ using a longest edge bisection algorithm. Set $i \leftarrow i + 1$ and go to step $\boldsymbol{1}$.
\end{algorithm}

\subsection{\FF{Numerical Experiments}}

We now provide two numerical experiments. In both examples, we consider problems where an exact solution can be obtained: succinctly, we fix the optimal state and adjoint state variables and compute the exact optimal control, the desired state $y_\Omega$, and the source term $f$. 

~\\
\textbf{Example 1 (L-shaped domain).} We set $\Omega=(-1,1)^2\setminus[0,1)\times(-1,0]$, $\texttt{a}=0.01$, $\texttt{b}=5$, and $\alpha=0.1$. The exact optimal state and adjoint state are given, in polar coordinates $(\rho,\omega)$ with $\omega\in[0,3\pi/2]$, by
\begin{align*}
\bar{y}&=3\sin(\pi(\rho\sin(\omega)+1)/2)\sin(\pi(\rho \cos(\omega)+1)/2)\rho^{2/3}\sin(2\omega/3), 
\\
\bar{p}&=2\cos(\pi\rho\sin(\omega)/2)\sin(\pi(\rho\cos(\omega)+1)/2)\rho^{2/3}\sin(2\omega/3).
\end{align*}

The purpose of this example is to investigate the performance of the devised a posteriori error estimators in a non--convex domain. 

\begin{figure}[!ht]
\centering
\psfrag{ndof}{$\mathsf{Ndof}^{-1/2}$}
\psfrag{ndof3}{\FF{$\mathsf{Ndof}^{-1/3}$}}
\psfrag{ndof-05}{$\mathsf{Ndof}^{-1/2}$}
\psfrag{ndof-1}{$\mathsf{Ndof}^{-1}$}
\psfrag{error}{$\|e\|_{\Omega}$}
\psfrag{estimator}{$\mathcal{E}_{ocp}$}
\psfrag{estimator-v}{$\mathsf{E}_{ocp}$}
\psfrag{error-y-unif}{\FF{Fully (unif)}}
\psfrag{error-y-unif-v}{\FF{Semi (unif)}}
\psfrag{error-y-adap}{\FF{Fully (adap)}}
\psfrag{error-y-var}{\FF{Semi (adap)}}
\psfrag{error-p-unif}{\FF{Fully (unif)}}
\psfrag{error-p-unif-v}{\FF{Semi (unif)}}
\psfrag{error-p-adap}{\FF{Fully (adap)}}
\psfrag{error-p-var}{\FF{Semi (adap)}}
\psfrag{error-u-unif}{\FF{Fully (unif)}}
\psfrag{error-u-unif-v}{\FF{Semi (unif)}}
\psfrag{error-u-adap}{\FF{Fully (adap)}}
\psfrag{error-u-var}{\FF{Semi (adap)}}
\psfrag{control-error}{$\|e_{\bar{u}}\|_{L^2(\Omega)}$}
\psfrag{eff-in}{{$\mathcal{I}_{eff}$}}
\psfrag{eff-in-v}{{$\mathsf{I}_{eff}$}}
\begin{minipage}[b]{0.32\textwidth}\centering
\scriptsize{\qquad \quad$\|\nabla e_{\bar{y}}\|_{L^2(\Omega)}$}\\
\includegraphics[trim={0 0 0 0},clip,width=3.9cm,height=3.6cm,scale=0.45]{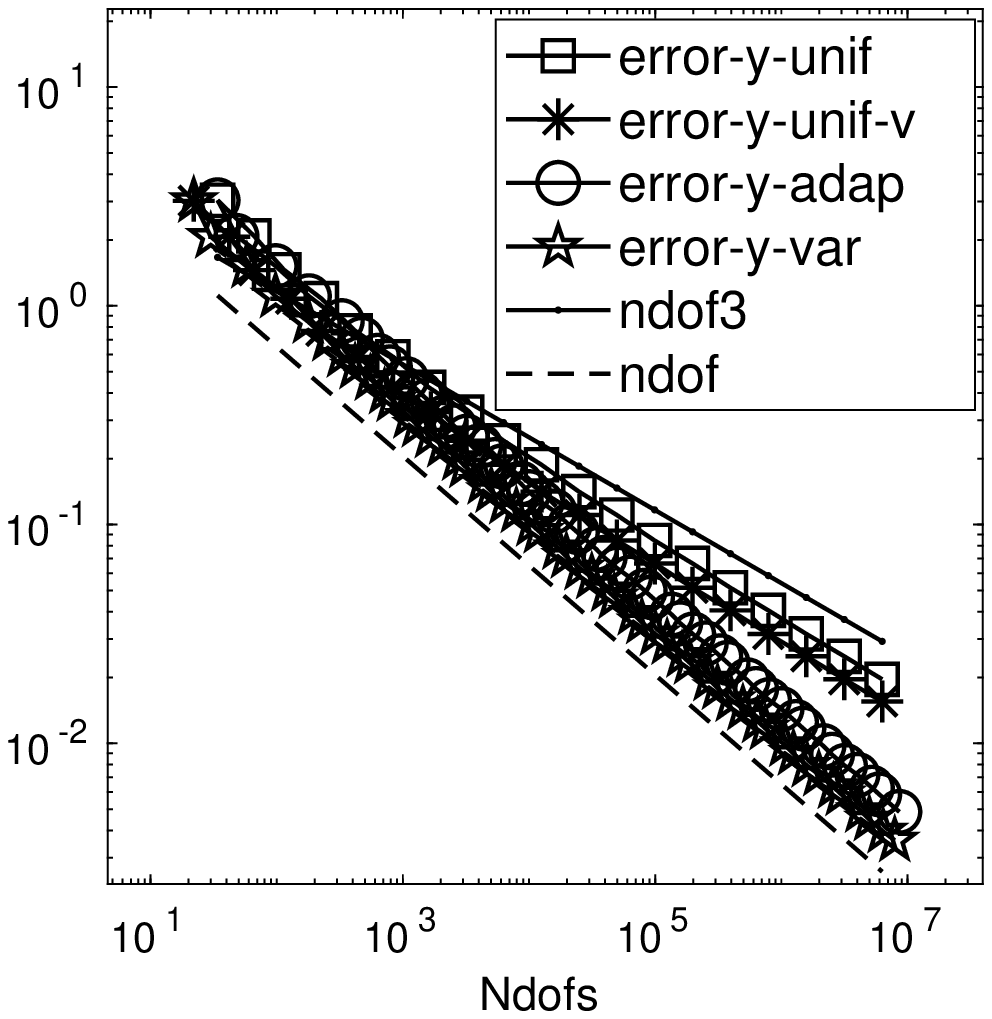} \\
\tiny{(A.1)}
\end{minipage}
\begin{minipage}[b]{0.32\textwidth}\centering
\scriptsize{\qquad \quad$\|\nabla e_{\bar{p}}\|_{L^2(\Omega)}$}\\
\includegraphics[trim={0 0 0 0},clip,width=3.9cm,height=3.6cm,scale=0.45]{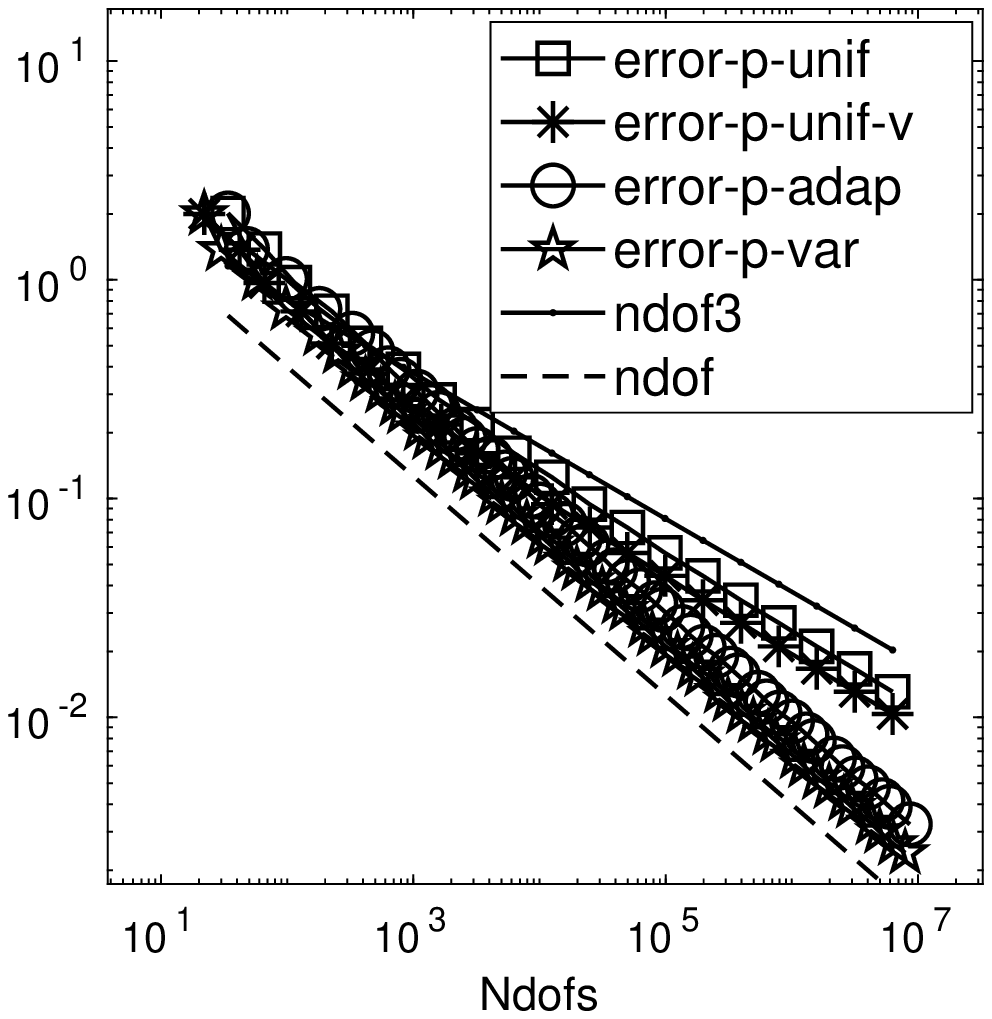} \\
\qquad \tiny{(A.2)}
\end{minipage}
\begin{minipage}[b]{0.32\textwidth}\centering
\scriptsize{\qquad$\|e_{\bar{u}}\|_{L^2(\Omega)}$ and $\|\mathsf{e}_{\bar{u}}\|_{L^2(\Omega)}$}\\
\includegraphics[trim={0 0 0 0},clip,width=3.9cm,height=3.6cm,scale=0.45]{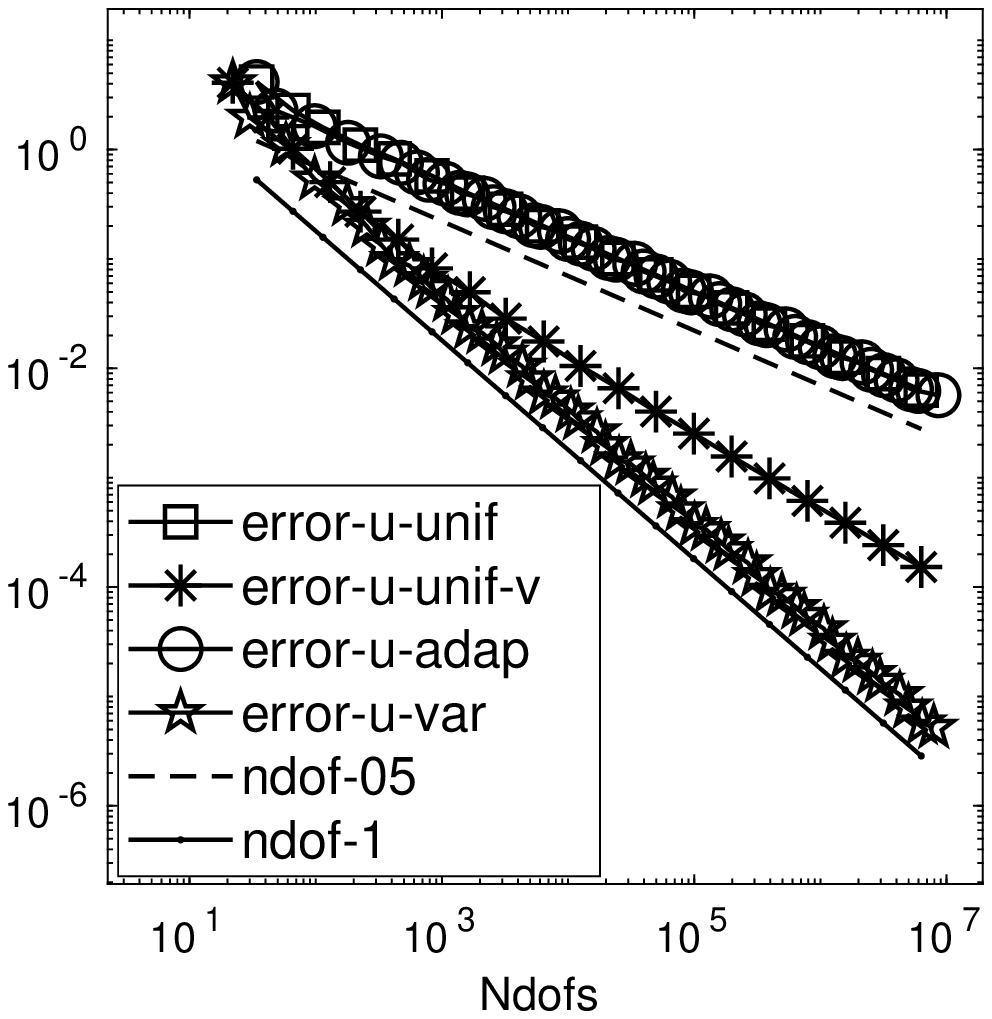} \\
\qquad \tiny{(A.3)}
\end{minipage}
\\~\\
\begin{minipage}[b]{0.32\textwidth}\centering
\scriptsize{\qquad Error vs Estimator (Fully)}\\
\includegraphics[trim={0 0 0 0},clip,width=3.9cm,height=3.6cm,scale=0.45]{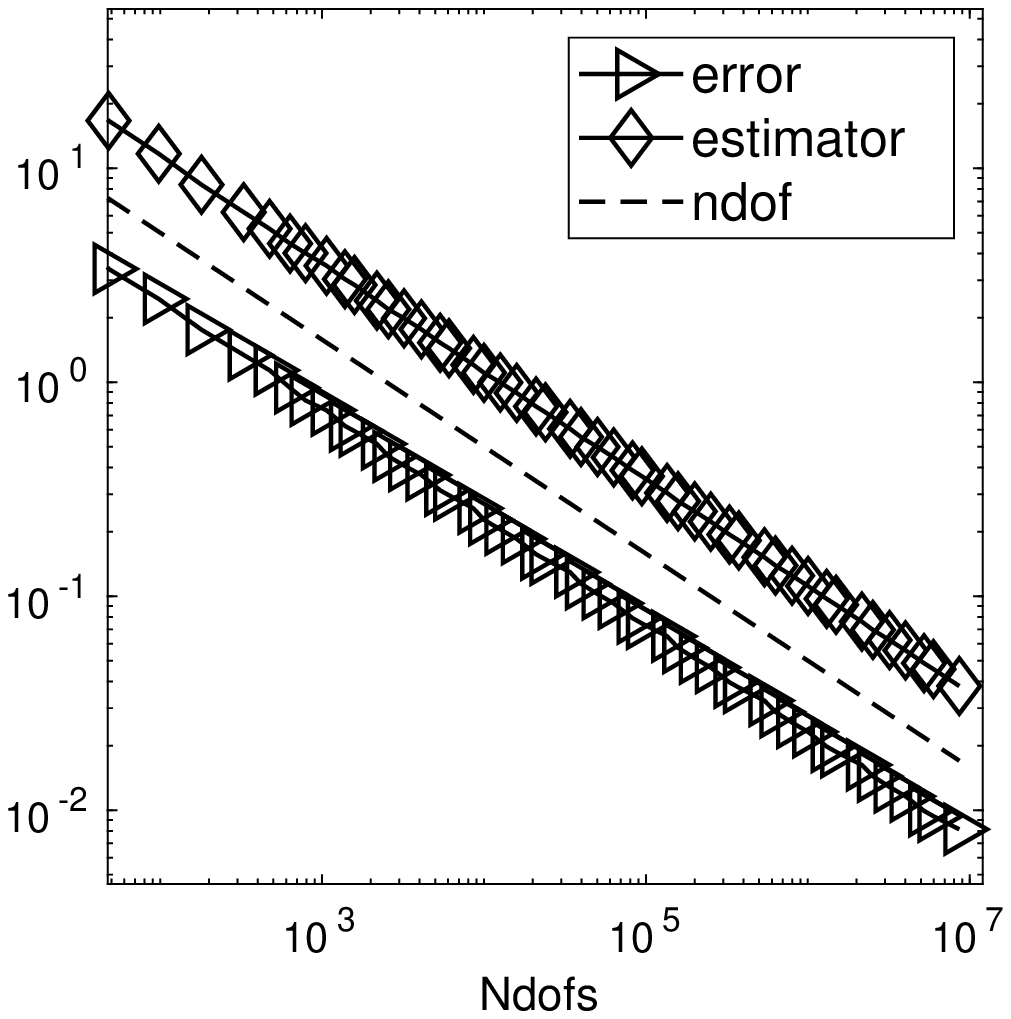}\\
\qquad \tiny{(A.4)}
\end{minipage}
\begin{minipage}[b]{0.32\textwidth}\centering
\scriptsize{\qquad Error vs Estimator (Semi)}\\
\psfrag{error}{$\|\mathsf{e}\|_{\Omega}$}
\includegraphics[trim={0 0 0 0},clip,width=3.9cm,height=3.6cm,scale=0.45]{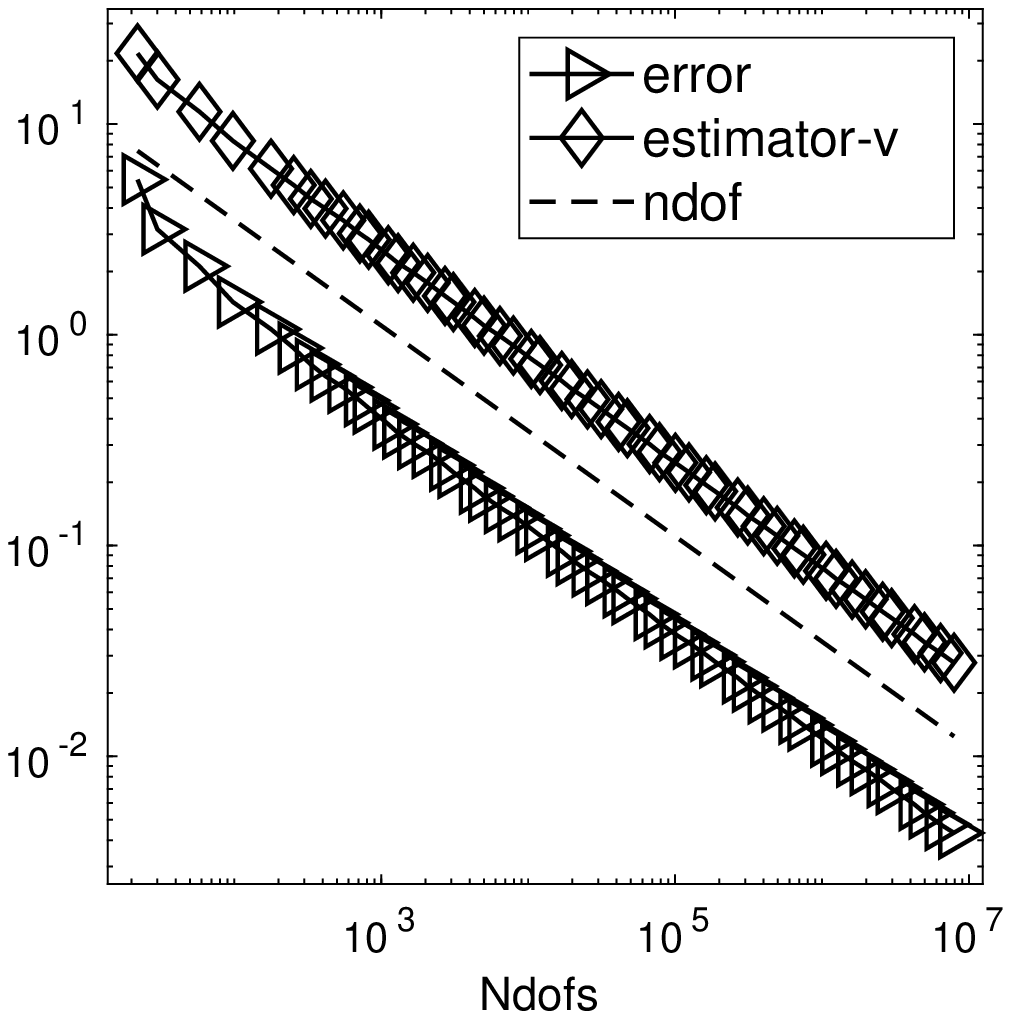}\\
\qquad \tiny{(A.5)}
\end{minipage}
\begin{minipage}[b]{0.32\textwidth}\centering
\scriptsize{Effectivity index}\\
\includegraphics[trim={0 0 0 0},clip,width=4.3cm,height=3.6cm,scale=0.45]{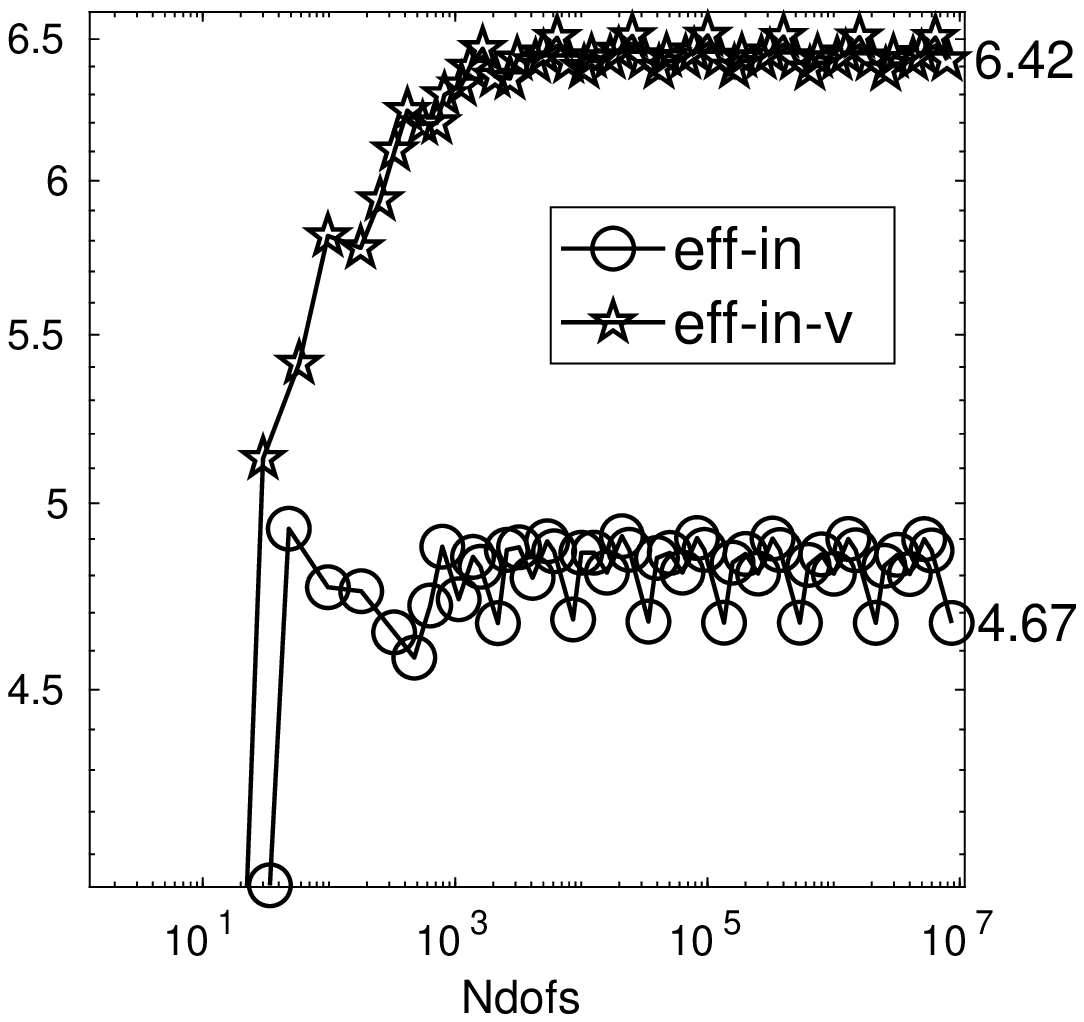}\\
\qquad \tiny{(A.6)}
\end{minipage}
\caption{Example 1: Experimental rates of convergence for the individual errors $\|\nabla e_{\bar{y}}\|_{L^{2}(\Omega)},\|\nabla e_{\bar{p}}\|_{L^{2}(\Omega)}$, $\|e_{\bar{u}}\|_{L^{2}(\Omega)}$, and $\|\mathsf{e}_{\bar{u}}\|_{L^{2}(\Omega)}$, for uniform and adaptive refinement (A.1)--(A.3), the total errors $\|e\|_{\Omega}$ and $\|\mathsf{e}\|_{\Omega}$ and the error estimators $\mathcal{E}_{ocp,\T}$ and $\mathsf{E}_{ocp,\T}$ for adaptive refinement (A.4) and (A.5), and the effectivity indices (A.6).}
\label{fig:ex_1}
\end{figure}

In Fig. \ref{fig:ex_1} we present the results obtained for Example 1. We present, in subfigures (A.1)--(A.3), experimental rates of convergence for all the individual contributions of the total errors $\|e\|_\Omega$ and $\|\mathsf{e}\|_{\Omega}$ when uniform and adaptive refinements are considered \FF{within both discretization schemes, i.e., the fully discrete scheme (Fully) and the semi-discrete scheme (Semi)}. We observe that our adaptive loops \emph{outperform} uniform refinement. In addition, \EO{we observe that our adaptive loops yield optimal experimental rates of convergence for all the individual contributions of the total errors $\|e\|_\Omega$ (Fully) and $\|\mathsf{e}\|_{\Omega}$ (Semi).} We also observe, in subfigures (A.4) and (A.5), that the error estimators \EO{$\mathcal{E}_{ocp}$ (Fully) and $\mathsf{E}_{ocp}$ (Semi)} exhibit optimal rates of convergence 
% \FF{when the corresponding adaptive refinement is considered}. 
Moreover, in subfigure (A.6), it can be observed that when \FF{$\mathsf{Ndof}$} increases, the effectivity indices \EO{$\mathcal{I}_{eff}$ (Fully) and $\mathsf{I}_{eff}$ (Semi)} are stabilized around the values 4.7 and 6.4, respectively.

%%%%%%%%%%%%%%%%%%%%%%%%%%%%%%%%%%%%%%%%%%%%%%%%%%%%%%%
%%%%%%%%%%%%%%%%%     FIGURE 1   %%%%%%%%%%%%%%%%%%%%%%
%%%%%%%%%%%%%%%%%%%%%%%%%%%%%%%%%%%%%%%%%%%%%%%%%%%%%%%

~\\
\textbf{Example 2 (Convex domain).} We consider $\Omega=(0,1)^3$, $\alpha=10^{-3}$, $\texttt{a}=0.2$, and $\texttt{b}=3$. The exact optimal state and adjoint state are given by $
\bar{y}= x_{1}x_{2}x_{3}(1-x_{1})(1-x_{2})(1-x_{3})$ and 
$\bar{p}=\bar{y}\cdot\arctan(100(x_{1}-0.5))$, respectively.

\begin{figure}[!ht]
\centering
\psfrag{ndof}{$\mathsf{Ndof}^{-1/3}$}
\psfrag{ndof-05}{$\mathsf{Ndof}^{-1/3}$}
\psfrag{ndof-1}{$\mathsf{Ndof}^{-2/3}$}
\psfrag{error}{$\|e\|_{\Omega}$}
\psfrag{estimator}{$\mathcal{E}_{ocp}$}
\psfrag{estimator-v}{$\mathsf{E}_{ocp}$}
\psfrag{error-y-unif}{\FF{Fully (unif)}}
\psfrag{error-y-unif-v}{\FF{Semi (unif)}}
\psfrag{error-y-adap}{\FF{Fully (adap)}}
\psfrag{error-y-var}{\FF{Semi (adap)}}
\psfrag{error-p-unif}{\FF{Fully (unif)}}
\psfrag{error-p-unif-v}{\FF{Semi (unif)}}
\psfrag{error-p-adap}{\FF{Fully (adap)}}
\psfrag{error-p-var}{\FF{Semi (adap)}}
\psfrag{error-u-unif}{\FF{Fully (unif)}}
\psfrag{error-u-unif-v}{\FF{Semi (unif)}}
\psfrag{error-u-adap}{\FF{Fully (adap)}}
\psfrag{error-u-var}{\FF{Semi (adap)}}
\psfrag{control-error}{$\|e_{\bar{u}}\|_{L^2(\Omega)}$}
\psfrag{eff-in}{{$\mathcal{I}_{eff}$}}
\psfrag{eff-in-v}{{$\mathsf{I}_{eff}$}}
\begin{minipage}[b]{0.32\textwidth}\centering
\scriptsize{\qquad \quad$\|\nabla e_{\bar{y}}\|_{L^2(\Omega)}$}\\
\includegraphics[trim={0 0 0 0},clip,width=3.9cm,height=3.6cm,scale=0.45]{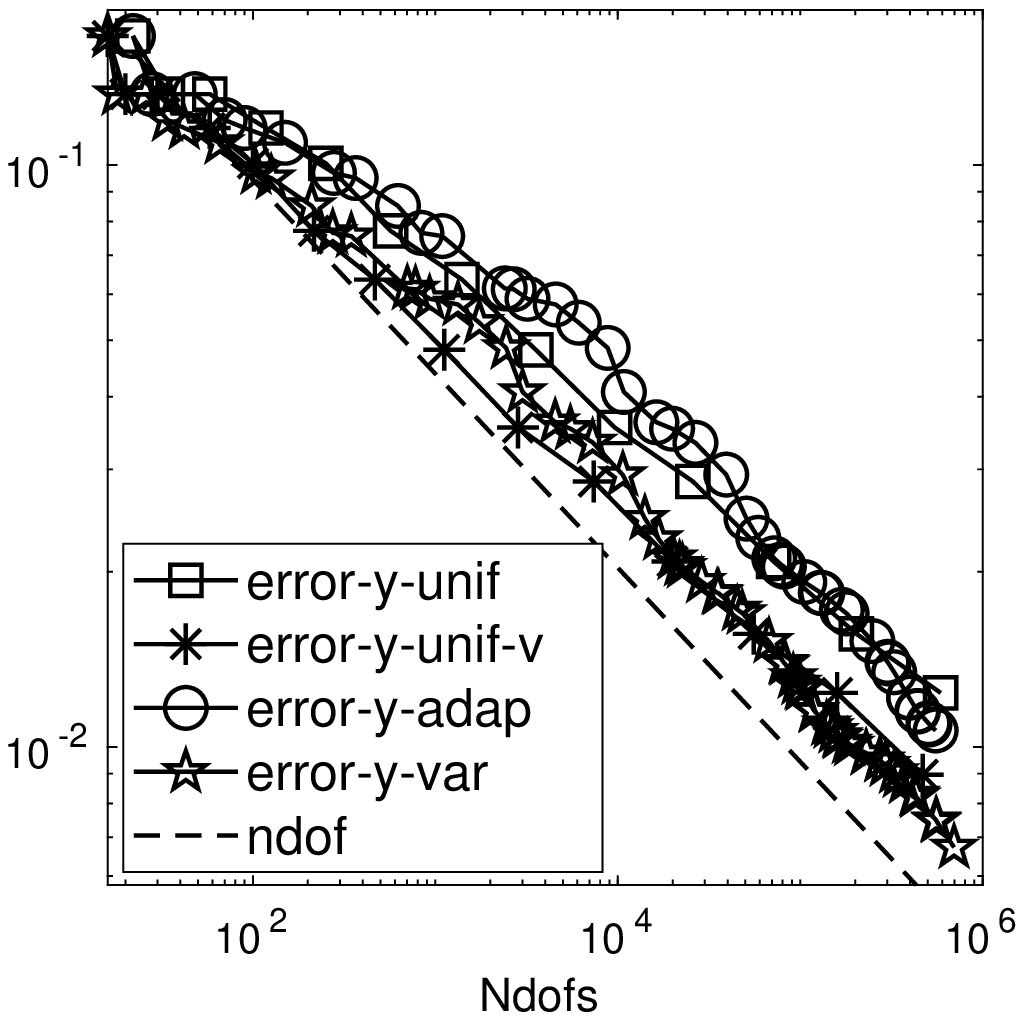} \\
\tiny{(B.1)}
\end{minipage}
\begin{minipage}[b]{0.32\textwidth}\centering
\scriptsize{\qquad \quad$\|\nabla e_{\bar{p}}\|_{L^2(\Omega)}$}\\
\includegraphics[trim={0 0 0 0},clip,width=3.9cm,height=3.6cm,scale=0.45]{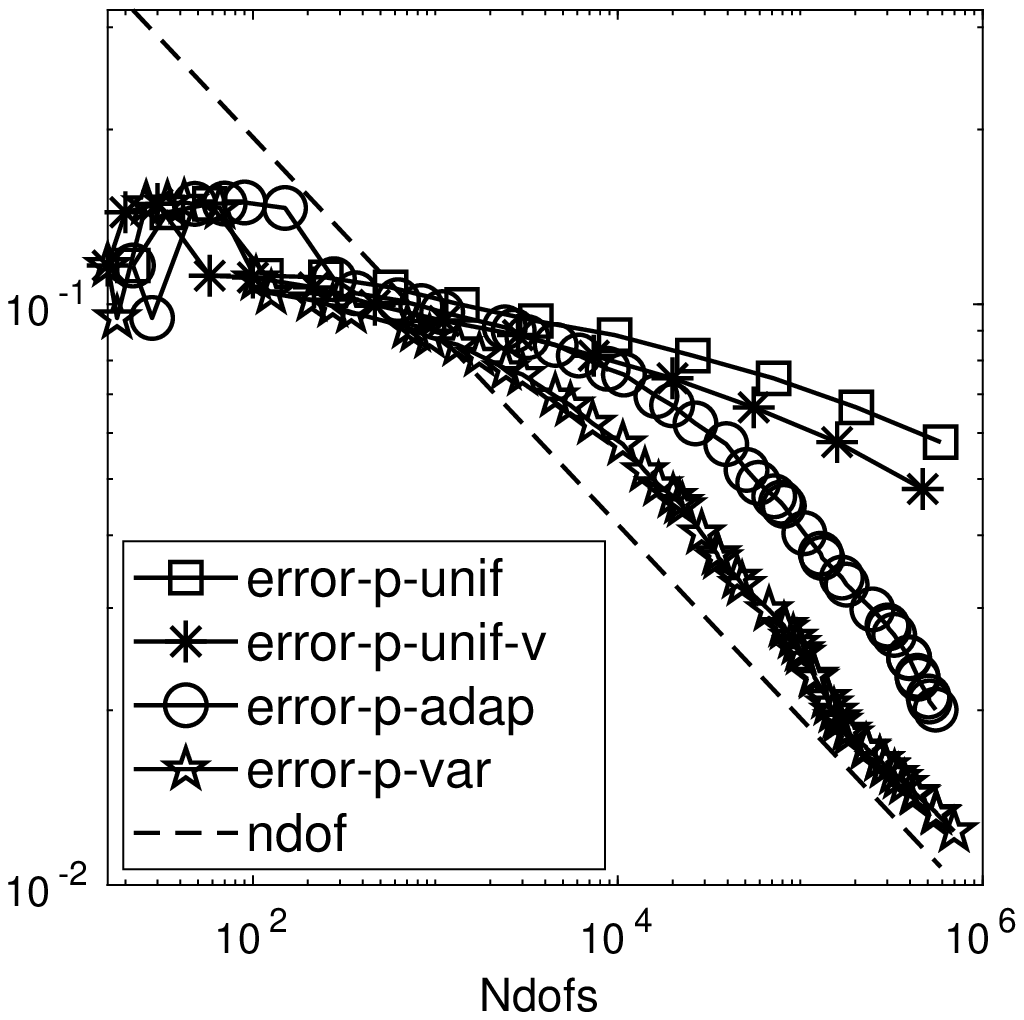} \\
\qquad \tiny{(B.2)}
\end{minipage}
\begin{minipage}[b]{0.32\textwidth}\centering
\scriptsize{\qquad$\|e_{\bar{u}}\|_{L^2(\Omega)}$ and $\|\mathsf{e}_{\bar{u}}\|_{L^2(\Omega)}$}\\
\includegraphics[trim={0 0 0 0},clip,width=3.9cm,height=3.6cm,scale=0.45]{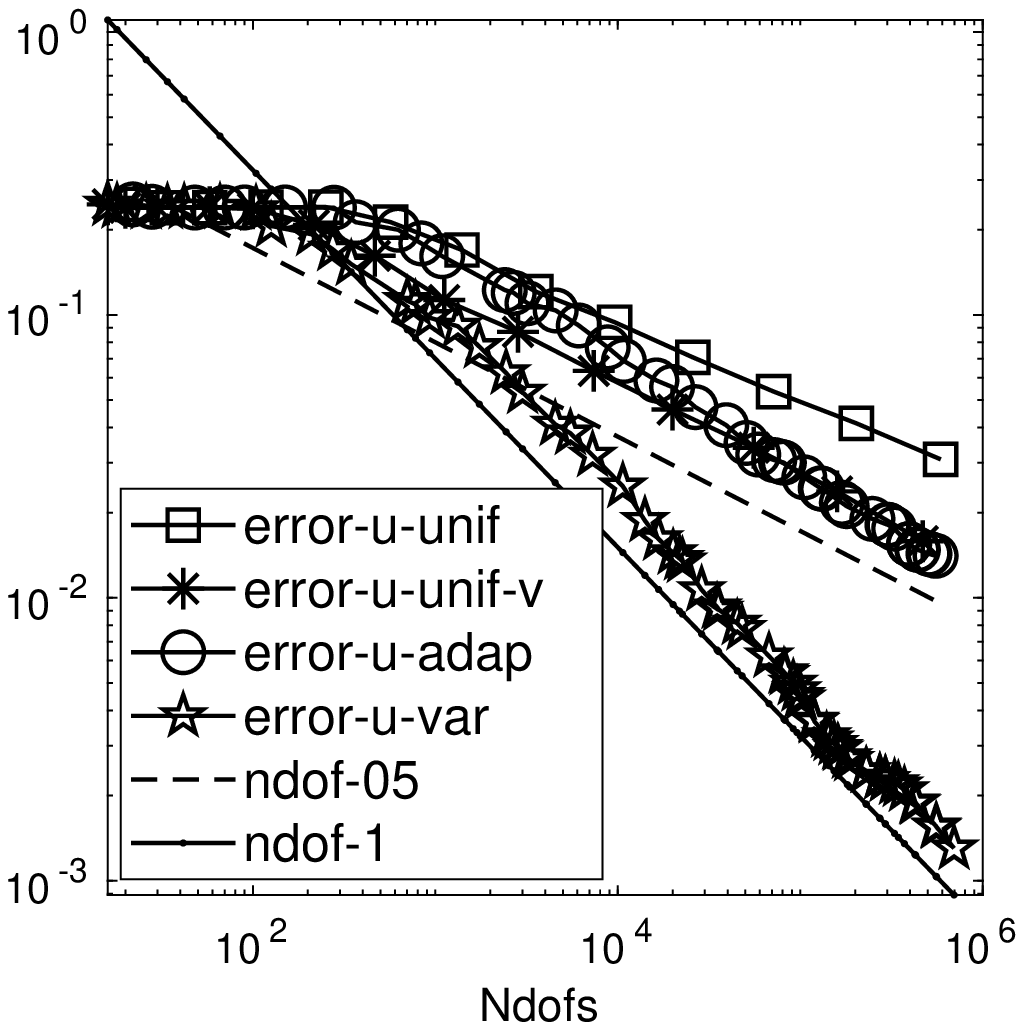} \\
\qquad \tiny{(B.3)}
\end{minipage}
\\~\\
\begin{minipage}[b]{0.32\textwidth}\centering
\scriptsize{\qquad Error vs Estimator (Fully)}\\
\includegraphics[trim={0 0 0 0},clip,width=3.9cm,height=3.6cm,scale=0.45]{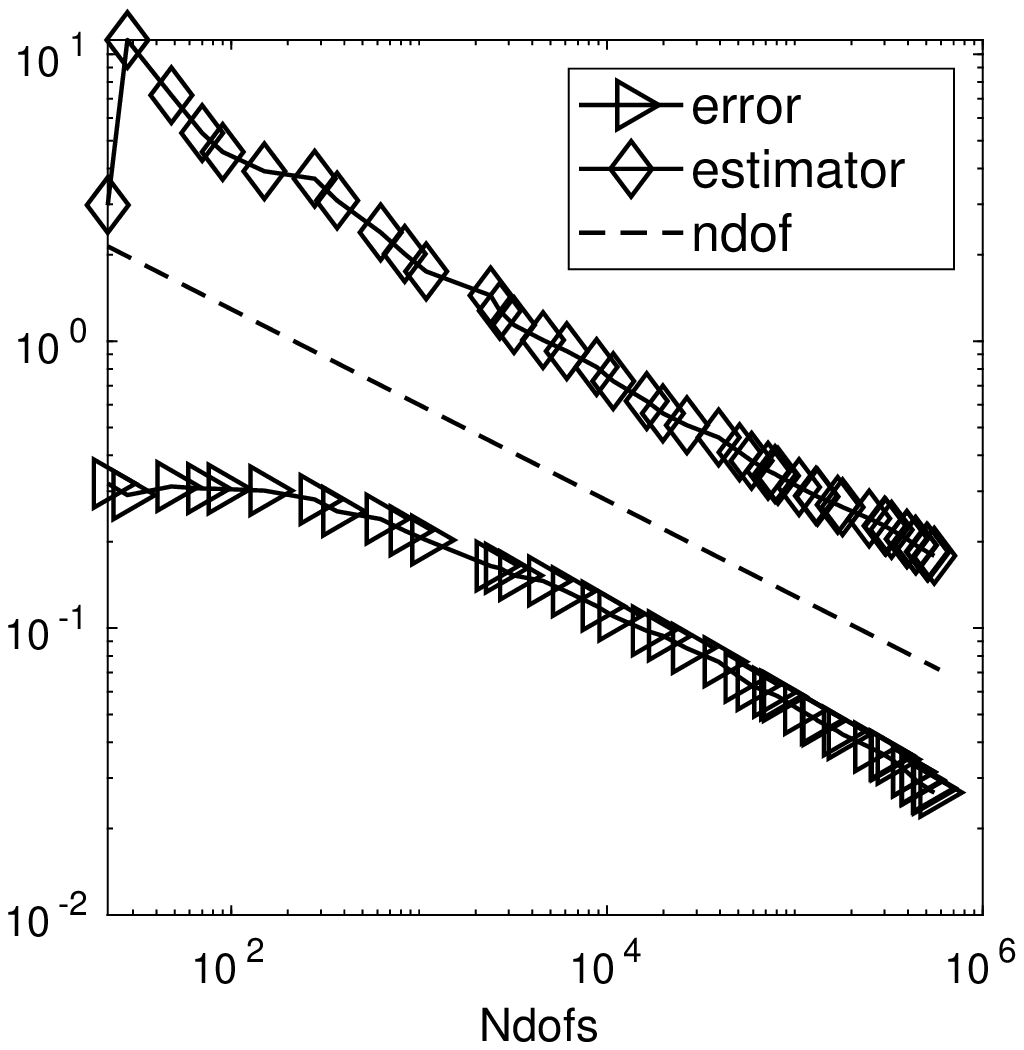}\\
\qquad \tiny{(B.4)}
\end{minipage}
\begin{minipage}[b]{0.32\textwidth}\centering
\scriptsize{\qquad Error vs Estimator (Semi)}\\
\psfrag{error}{$\|\mathsf{e}\|_{\Omega}$}
\includegraphics[trim={0 0 0 0},clip,width=3.9cm,height=3.6cm,scale=0.45]{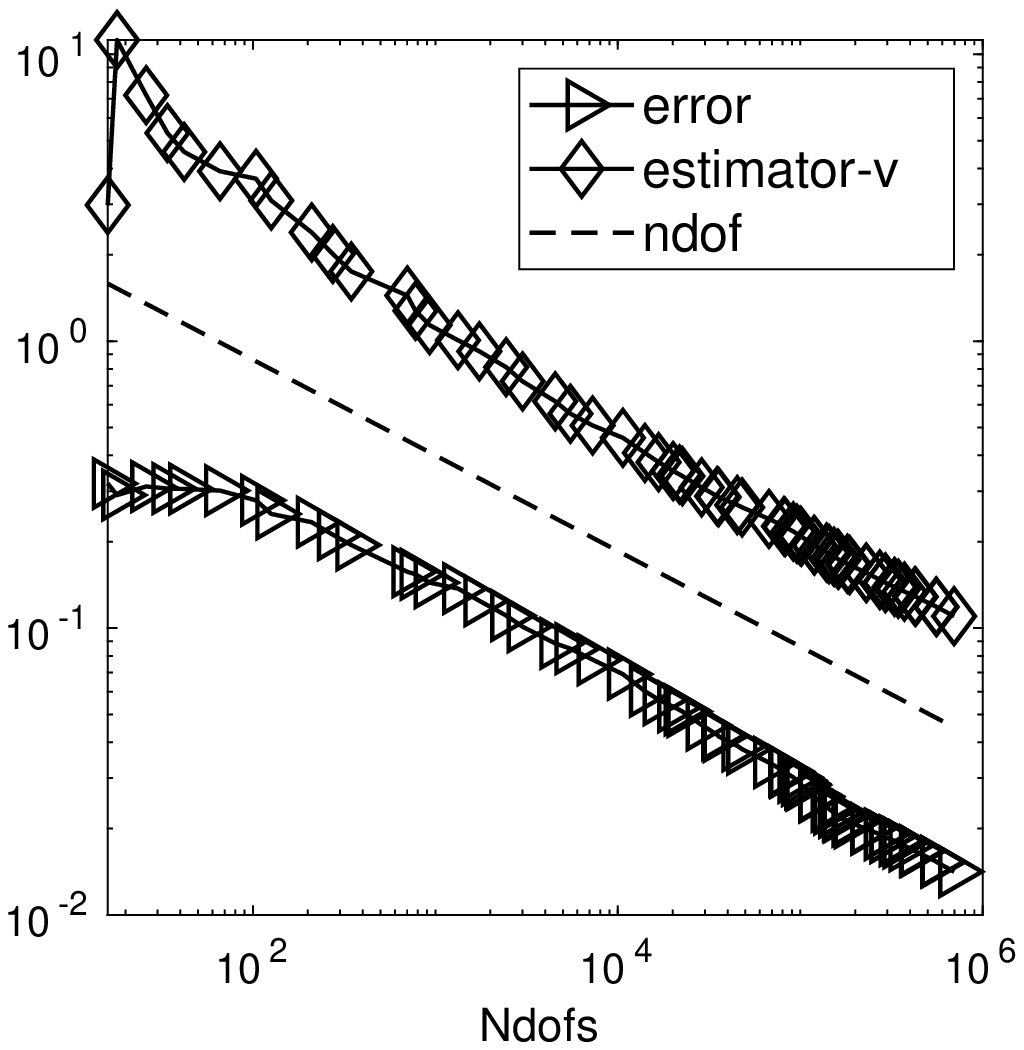}\\
\qquad \tiny{(B.5)}
\end{minipage}
\begin{minipage}[b]{0.32\textwidth}\centering
\scriptsize{Effectivity index}\\
\includegraphics[trim={0 0 0 0},clip,width=4.3cm,height=3.6cm,scale=0.45]{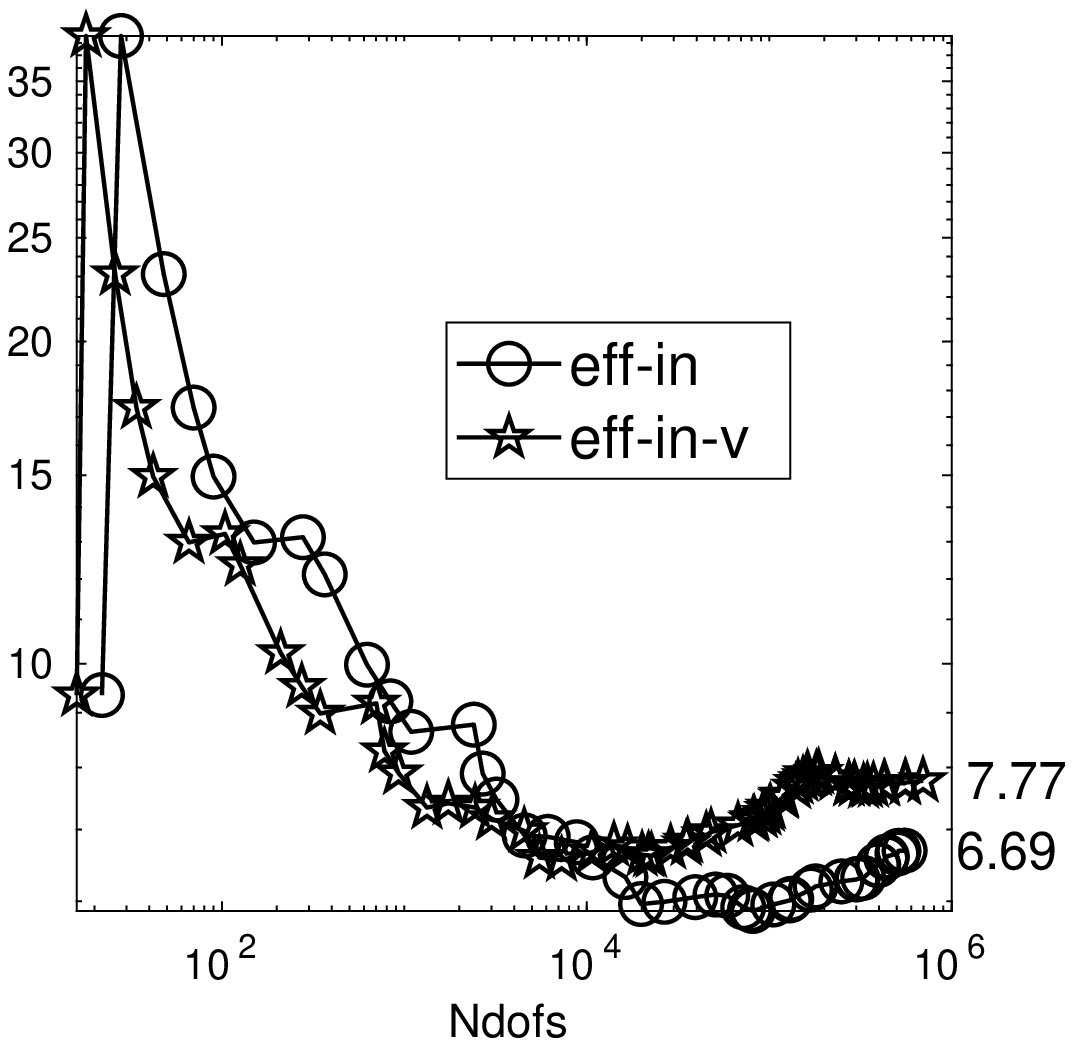}\\
\qquad \tiny{(B.6)}
\end{minipage}
\caption{Example 2: Experimental rates of convergence for the individual errors $\|\nabla e_{\bar{y}}\|_{L^{2}(\Omega)},\|\nabla e_{\bar{p}}\|_{L^{2}(\Omega)}$, $\|e_{\bar{u}}\|_{L^{2}(\Omega)}$, and $\|\mathsf{e}_{\bar{u}}\|_{L^{2}(\Omega)}$, for uniform and adaptive refinement (B.1)--(B.3), the total errors $\|e\|_{\Omega}$ and $\|\mathsf{e}\|_{\Omega}$ and error estimators $\mathcal{E}_{ocp,\T}$ and $\mathsf{E}_{ocp,\T}$ for adaptive refinement (B.4) and (B.5), and the effectivity indices (B.6).}
\label{fig:ex_2}
\end{figure}

In Fig. \ref{fig:ex_2} we present the results obtained for Example 2. Similar conclusions to the ones presented for Example 1 can be derived. In particular, we observe optimal experimental rates of convergence for all the individual contributions of the total errors \EO{$\|e\|_\Omega$ (Fully) and $\|\mathsf{e}\|_{\Omega}$ (Semi)} when adaptive refinement is considered.

%%%%%%%%%%%%%%%%%%%%%%%%%%%%%%%%%%%%%%%%%%%%%%%%%%%%%%%
%%%%%%%%%%%%%%%%%     FIGURE 2    %%%%%%%%%%%%%%%%%%%%%
%%%%%%%%%%%%%%%%%%%%%%%%%%%%%%%%%%%%%%%%%%%%%%%%%%%%%%%

\section{Conclusions}\label{sec:conclusions}

In this work, we designed AFEMs for an optimal control problem with a bilinear state equation; the control variable enters the state equation as a coefficient. Two different discretization schemes were considered to approximate an optimal solution: a fully discrete scheme and a semi-discrete one based on the so-called variation discretization approach. We constructed, for each one of these schemes, residual--type a posteriori error estimators that \FF{are formed by} the sum of contributions related to the discretization of the state and adjoint equations and, additionally, the discretization of the control variable for when the fully discrete scheme is considered. We derived global reliability and efficiency estimates for the error estimator associated to the fully discrete scheme whereas global reliability and \emph{local} efficiency estimates were derived for the error estimator associated to the semi-discrete scheme. Finally, we \EO{performed} numerical tests that showed competitive performances of both error estimators when used to drive adaptive procedures.

\begin{acknowledgements}
The first author was supported by UTFSM through Beca de Mantenci\'on. The second author was partially supported by ANID through FONDECYT Project \EO{1220156}.
\end{acknowledgements}

\footnotesize
\bibliographystyle{acm}

\bibliography{biblio_bilineal}

\begin{thebibliography}{10}

\bibitem{MR2424078}
{\sc Adams, R.~A., and Fournier, J. J.~F.}
\newblock {\em Sobolev spaces}, second~ed., vol.~140 of {\em Pure and Applied
  Mathematics (Amsterdam)}.
\newblock Elsevier/Academic Press, Amsterdam, 2003.

\bibitem{MR4122501}
{\sc Allendes, A., Fuica, F., and Ot\'{a}rola, E.}
\newblock Adaptive finite element methods for sparse {PDE}-constrained
  optimization.
\newblock {\em IMA J. Numer. Anal. 40}, 3 (2020), 2106--2142.

\bibitem{MUMPS1}
{\sc Amestoy, P., Duff, I., and L'Excellent, J.-Y.}
\newblock Multifrontal parallel distributed symmetric and unsymmetric solvers.
\newblock {\em Computer Methods in Applied Mechanics and Engineering 184}, 2
  (2000), 501 -- 520.

\bibitem{MUMPS2}
{\sc Amestoy, P.~R., Duff, I.~S., L'Excellent, J.-Y., and Koster, J.}
\newblock A fully asynchronous multifrontal solver using distributed dynamic
  scheduling.
\newblock {\em SIAM J. Matrix Anal. Appl. 23}, 1 (2001), 15--41 (electronic).

\bibitem{MR1780911}
{\sc Becker, R., Kapp, H., and Rannacher, R.}
\newblock Adaptive finite element methods for optimal control of partial
  differential equations: basic concept.
\newblock {\em SIAM J. Control Optim. 39}, 1 (2000), 113--132.

\bibitem{MR2373954}
{\sc Brenner, S.~C., and Scott, L.~R.}
\newblock {\em The mathematical theory of finite element methods}, third~ed.,
  vol.~15 of {\em Texts in Applied Mathematics}.
\newblock Springer, New York, 2008.

\bibitem{MR3174031}
{\sc Chang, Y., Yang, D., and Zhang, Z.}
\newblock Adaptive finite element approximation for a class of parameter
  estimation problems.
\newblock {\em Appl. Math. Comput. 231\/} (2014), 284--298.

\bibitem{MR3103238}
{\sc Chen, Y., Lu, Z., and Huang, Y.}
\newblock Superconvergence of triangular {R}aviart-{T}homas mixed finite
  element methods for a bilinear constrained optimal control problem.
\newblock {\em Comput. Math. Appl. 66}, 8 (2013), 1498--1513.

\bibitem{MR0520174}
{\sc Ciarlet, P.~G.}
\newblock {\em The finite element method for elliptic problems}.
\newblock North-Holland Publishing Co., Amsterdam-New York-Oxford, 1978.
\newblock Studies in Mathematics and its Applications, Vol. 4.

\bibitem{MR2971171}
{\sc Clason, C., and Jin, B.}
\newblock A semismooth {N}ewton method for nonlinear parameter identification
  problems with impulsive noise.
\newblock {\em SIAM J. Imaging Sci. 5}, 2 (2012), 505--538.

\bibitem{MR3693332}
{\sc Fu, H., Guo, H., Hou, J., and Zhang, J.}
\newblock A stabilized mixed finite element approximation of bilinear state
  optimal control problems.
\newblock {\em Comput. Math. Appl. 74}, 6 (2017), 1246--1261.

\bibitem{MR3621827}
{\sc Gong, W., and Yan, N.}
\newblock Adaptive finite element method for elliptic optimal control problems:
  convergence and optimality.
\newblock {\em Numer. Math. 135}, 4 (2017), 1121--1170.

\bibitem{MR2434065}
{\sc Hinterm\"{u}ller, M., Hoppe, R. H.~W., Iliash, Y., and Kieweg, M.}
\newblock An a posteriori error analysis of adaptive finite element methods for
  distributed elliptic control problems with control constraints.
\newblock {\em ESAIM Control Optim. Calc. Var. 14}, 3 (2008), 540--560.

\bibitem{MR2122182}
{\sc Hinze, M.}
\newblock A variational discretization concept in control constrained
  optimization: the linear-quadratic case.
\newblock {\em Comput. Optim. Appl. 30}, 1 (2005), 45--61.

\bibitem{MR2516528}
{\sc Hinze, M., Pinnau, R., Ulbrich, M., and Ulbrich, S.}
\newblock {\em Optimization with {PDE} constraints}, vol.~23 of {\em
  Mathematical Modelling: Theory and Applications}.
\newblock Springer, New York, 2009.

\bibitem{MR2843956}
{\sc Hinze, M., and Tr\"{o}ltzsch, F.}
\newblock Discrete concepts versus error analysis in {PDE}-constrained
  optimization.
\newblock {\em GAMM-Mitt. 33}, 2 (2010), 148--162.

\bibitem{MR3212590}
{\sc Kohls, K., R\"{o}sch, A., and Siebert, K.~G.}
\newblock A posteriori error analysis of optimal control problems with control
  constraints.
\newblock {\em SIAM J. Control Optim. 52}, 3 (2014), 1832--1861.

\bibitem{MR2536007}
{\sc Kr\"{o}ner, A., and Vexler, B.}
\newblock A priori error estimates for elliptic optimal control problems with a
  bilinear state equation.
\newblock {\em J. Comput. Appl. Math. 230}, 2 (2009), 781--802.

\bibitem{MR2680928}
{\sc Kunisch, K., Liu, W., Chang, Y., Yan, N., and Li, R.}
\newblock Adaptive finite element approximation for a class of parameter
  estimation problems.
\newblock {\em J. Comput. Math. 28}, 5 (2010), 645--675.

\bibitem{MR1887737}
{\sc Liu, W., and Yan, N.}
\newblock A posteriori error estimates for distributed convex optimal control
  problems.
\newblock {\em Adv. Comput. Math. 15}, 1-4 (2001), 285--309.

\bibitem{Nochetto_etal2009}
{\sc Nochetto, R.~H., Siebert, K.~G., and Veeser, A.}
\newblock Theory of adaptive finite element methods: an introduction.
\newblock In {\em Multiscale, nonlinear and adaptive approximation}. Springer,
  Berlin, 2009, pp.~409--542.

\bibitem{Troltzsch}
{\sc Tr\"oltzsch, F.}
\newblock {\em Optimal control of partial differential equations}, vol.~112 of
  {\em Graduate Studies in Mathematics}.
\newblock American Mathematical Society, Providence, RI, 2010.
\newblock Theory, methods and applications, Translated from the 2005 German
  original by J\"urgen Sprekels.

\bibitem{MR3059294}
{\sc Verf\"{u}rth, R.}
\newblock {\em A posteriori error estimation techniques for finite element
  methods}.
\newblock Numerical Mathematics and Scientific Computation. Oxford University
  Press, Oxford, 2013.

\bibitem{MR2373479}
{\sc Vexler, B., and Wollner, W.}
\newblock Adaptive finite elements for elliptic optimization problems with
  control constraints.
\newblock {\em SIAM J. Control Optim. 47}, 1 (2008), 509--534.

\end{thebibliography}
\end{document}